\documentclass[a4paper,11pt]{article}
\usepackage[utf8]{inputenc}
\usepackage[english]{babel}
\usepackage[T1]{fontenc}
\usepackage{libertine}
\usepackage{graphicx}
\usepackage{comment}
\usepackage{amsthm,amsmath,amsfonts,stmaryrd,mathrsfs,amssymb}
\usepackage{mathtools}
\usepackage{bbm}
\usepackage{enumitem}
\usepackage{chemarrow}
\usepackage{MnSymbol}
\usepackage{tikz}
\usepackage[top=2.5cm, bottom=2.5cm, left=2cm, right=2cm]{geometry}
\usepackage{subfig}
\usepackage{multirow}
\usepackage{physics}
\usepackage{array}
\usepackage{bigints}
\usepackage{etoolbox}
\usepackage{mdframed}
\usepackage{color}
\usepackage{accents}
\usepackage[unicode,colorlinks=true,linkcolor=blue,citecolor=red,pdfencoding=auto,psdextra]{hyperref}
\usepackage{tabulary}
\usepackage{refcount} 
\usepackage{mleftright}\mleftright 
\usepackage{booktabs}
\makeatletter
\newcommand{\mylabel}[2]{#2\def\@currentlabel{#2}\label{#1}}
\makeatother

\newlength{\dhatheight}
\newcommand{\doublehat}[1]{%
	\settoheight{\dhatheight}{\ensuremath{\hat{#1}}}%
	\addtolength{\dhatheight}{-0.35ex}%
	\hat{\vphantom{\rule{1pt}{\dhatheight}}%
		\smash{\hat{#1}}}}

\newcommand{\ensemblenombre}[1]{\mathbb{#1}}

\newcommand{\Z}{\ensemblenombre{Z}}

\newcommand{\R}{\ensemblenombre{R}}

\newcounter{stepgroup}
\newcounter{step}

\AtBeginEnvironment{proof}{\stepcounter{stepgroup}\setcounter{step}{0}}
\newcommand{\step}[1]{\leavevmode\refstepcounter{step}\emph{Step~\thestep: #1.}}
\tikzset{
 vin/.style={circle,fill=black,inner sep=1.7pt},
 vev/.style={circle,draw=black,fill=gray!45,inner sep=1.7pt},
 vbd/.style={circle,draw=black,fill=white,inner sep=1.7pt},
 vroot/.style={rectangle,fill=black,inner sep=2.2pt},
 vghost/.style={circle,draw=gray!60,fill=white,inner sep=1.7pt},
 openlong/.style={thick},
 closedlong/.style={densely dotted,thick},
 shortedge/.style={gray!55},
 shortdash/.style={dashed,thick},
 subtreestyle/.style={fill=gray!15,draw=gray!55},
}
\newcommand{\gfcoords}{%
 \coordinate (r0) at (0,-1); \coordinate (r) at (0,0);
 \coordinate (a) at (-1.9,1); \coordinate (b) at (0.7,1); \coordinate (c) at (2.1,1);
 \coordinate (a1) at (-2.8,2); \coordinate (a2) at (-1.1,2); \coordinate (b1) at (0.7,2);
 \coordinate (a11) at (-3.3,3); \coordinate (a12) at (-2.3,3);}
\newcommand{\subtreeat}[1]{\filldraw[subtreestyle] (#1) -- ++(-0.38,0.85) -- ++(0.76,0) -- cycle;}

\newcommand{\Ec}[1]{\mathbb{E} \mathinner{\left[#1\right]}}

\newcommand{\Pp}[1]{\mathbb{P} \mathinner{\left(#1\right)}}

\newcommand{\Ecsq}[2]{\mathbb{E} \mathinner{\left[#1\mathrel{}\middle|\mathrel{}#2\right]}}

\newcommand{\Ppsq}[2]{\mathbb{P} \mathinner{\left(#1\mathrel{}\middle|\mathrel{}#2\right)}}

\newcommand{\Var}[1]{\mathrm{Var}\mathinner{\left(#1\right)}}

\newcommand{\intervalle}[4]{\left#1#2
	\mathclose{}\mathpunct{},#3
	\right#4}
\newcommand{\intervalleff}[2]{\intervalle{[}{#1}{#2}{]}}
\newcommand{\intervalleof}[2]{\intervalle{(}{#1}{#2}{]}}
\newcommand{\intervallefo}[2]{\intervalle{[}{#1}{#2}{)}}
\newcommand{\intervalleoo}[2]{\intervalle{(}{#1}{#2}{)}}

\newcommand{\enstq}[2]{\mathinner{\left\lbrace#1\mathrel{}\middle|\mathrel{}#2\right\rbrace}}

\newcommand{\ind}[1]{\mathbf{1}_{\lbrace #1 \rbrace}}

\DeclareMathOperator{\haut}{ht}
\DeclareMathOperator{\supp}{supp}
\DeclareMathOperator{\GF}{GP} 

\newtheorem{theorem}{Theorem}
\newtheorem{proposition}[theorem]{Proposition}
\newtheorem{lemma}[theorem]{Lemma}

\newtheorem{corollary}[theorem]{Corollary}
\newtheorem{remark}[theorem]{Remark}

\begin{document}
	\title{Percolation on grandparent graphs of Bienaym\'e trees}
	\author{Laurent \textsc{M\'enard}\thanks{MODAL'X, UPL, Univ. Paris Nanterre, CNRS, F92000 Nanterre France.\hfill  \href{mailto:lmenard@parisnanterre.fr}{\texttt{lmenard@parisnanterre.fr}}}
		\quad and \quad Delphin \textsc{S\'enizergues}\thanks{MODAL'X, UPL, Univ. Paris Nanterre, CNRS, F92000 Nanterre France.\hfill  \href{mailto:dsenizer@parisnanterre.fr}{\texttt{dsenizer@parisnanterre.fr}}}}
	\maketitle

\begin{abstract}
We consider Bernoulli bond percolation on the grandparent graph of a Bienaym\'e tree, obtained by adding to the tree an edge between every vertex and its grandparent. 
We show that the critical parameter $p_c$ is the unique solution of an explicit equation involving the generating function of the offspring distribution, and we study how $p_c$ depends on this distribution beyond its mean $m$. 
When $m\geq 1.569$, the offspring distribution concentrated on the (one or two) integers closest to $m$ is the unique maximiser of $p_c$ among all offspring distributions with mean $m$, and its critical parameter is $\frac1{m^2}-\frac1{m^3}+O\left(\frac1{m^4}\right)$ as $m\to\infty$; 
more generally, 
for two distributions that are comparable in the convex order,
$p_c$ is larger for the smaller of the two distributions under the condition that its support bounded some explicit $N(m)$ of order $m^4$.
These restrictions cannot be removed: the convex order alone does not determine the order of the critical parameters, even among distributions with support bounded by $O(m^4)$, and the extremal property fails for $m$ close to $1$. 
Finally, for every $m>1$, the critical parameter $p_c$ can be arbitrarily close to $0$ for general offspring distributions with mean $m$.
If we only consider offspring distributions that do not charge $0$, so that the obtained tree is almost surely infinite, this remains true if and only if $m\geq 2$.
\end{abstract}

\section{Introduction}

\subsection{Model and questions}\label{sec:intro model}

Let $X$ be a random variable with values in $\Z_+=\{0,1,2,\dots\}$ and mean $m:=\Ec{X}\in\intervalleoo{1}{\infty}$, and let $T$ be a Bienaym\'e tree (also called Galton--Watson tree) whose offspring distribution is the law of $X$. The \emph{grandparent graph} $\GF(T)$ of $T$ is the graph obtained from $T$ by adding an edge between every vertex and its grandparent, when it exists (see Section~\ref{sec:determination} for a formal definition). We study Bernoulli bond percolation with parameter $p\in\intervalleff{0}{1}$ on $\GF(T)$ and define the critical parameter
\[
p_c(X):=\inf\enstq{p\in\intervalleff{0}{1}}{\Pp{\text{the open cluster of the root of } \GF(T) \text{ is infinite}}>0},
\]
where the probability is taken with respect to both the tree and the percolation (annealed setting). It only depends on the law of $X$. 
By a standard $0$--$1$ argument (see Appendix~\ref{sec:annealed quenched}), for almost every tree $T$ on the event of survival, the critical parameter of percolation on the fixed graph $\GF(T)$ is equal to $p_c(X)$: annealed and quenched critical parameters coincide.

Since $T$ is a subgraph of $\GF(T)$, we have $p_c(X)\leq 1/m$; using only the edges between vertices and their grandparent, one easily sees that $p_c(X)\leq 1/m^2$. 
For bond percolation on $T$ itself, Lyons~\cite{Lyo90} showed that the critical parameter is $1/m$: it only depends on the mean of the offspring distribution. 
On the contrary, $p_c(X)$ depends on the whole law of $X$, and the main question we address is the following: \emph{at fixed mean $m$, how does $p_c(X)$ depend on the law of $X$? Which laws maximise or minimise $p_c$?}

A natural way to compare distributions with the same mean is the convex order \cite{SS07}: for integrable $\Z_+$-valued random variables $X$ and $Y$, we write 
$X\preceq_{\mathrm{cx}}Y$ if $\Ec{f(X)}\leq\Ec{f(Y)}$ for every convex function $f:\Z_+\to\R$ such that both expectations are well defined. 
Applied to the functions $k\mapsto k$ and $k\mapsto -k$, this definition implies $\Ec{X}=\Ec{Y}$. 
The relation $X\preceq_{\mathrm{cx}}Y$ expresses the fact that $Y$ is ``more variable'' than $X$. 
We write $X\prec_{\mathrm{cx}}Y$ if $X\preceq_{\mathrm{cx}}Y$ and additionally $X$ and $Y$ do not have the same law. 

For $m>1$ we denote by $X^{(m)}$ a random variable with law
\begin{align}\label{eq:def loi X^{(m)}}
(1-\{m\})\cdot\delta_{\lfloor m\rfloor}+\{m\}\cdot\delta_{\lfloor m\rfloor+1},\qquad\text{where } \{m\}:=m-\lfloor m\rfloor,
\end{align}
which is the smallest element for the convex order among $\Z_+$-valued random variables with mean $m$. 
When $m$ is an integer, $X^{(m)}=m$ almost surely and $T$ is the $m$-ary tree.

\subsection{Main result: \texorpdfstring{$p_c$}{p\_c} decreases along the convex order}
Our main result allows us to compare the critical parameters $p_c(X)$ and $p_c(Y)$ when $X\mathrel{\prec_{\mathrm{cx}}} Y$, under some additional assumption on the support of $X$, when $m$ is not too close to $1$.
\begin{theorem}\label{thm:convex order + bound support + large enough mean = comparison}
For $m>1$, set
	\[
	N(m):=\left\lfloor 1+\frac{(m-1)^2(m^3+1)}{m}\right\rfloor .
	\]
	If $X$ and $Y$ are $\Z_+$-valued random variables with mean $m$ such that $X \mathrel{\prec_{\mathrm{cx}}} Y$ and $\supp X \subset \{0, \dots, N(m)\}$, then
	\begin{align*}
		p_c(X) > p_c(Y).
	\end{align*}
\end{theorem}
\begin{remark}\label{rem:N(m)}
Note that the statement is empty when $N(m)\leq 1$, which happens if and only if $m<m_1$, where $m_1\approx 1.5682$ is the unique solution in $\intervalleoo{1}{2}$ of $(m-1)^2(m^3+1)=m$.
Moreover, as we check in the proof of Corollary~\ref{cor:maximal percolation threseshold}, we have $N(m)\geq\lceil m\rceil$ for every $m\geq m_1$, ensuring that the statement is not empty in this case.
When $m$ is large, we have $N(m)=m^4-2m^3+O(m^2)$.
\end{remark}

The theorem above does not allow us to compare all possible pairs of offspring distributions comparable in the convex order, but it is enough to identify the distributions with the highest $p_c$, at least when $m$ is not too close to $1$.
Recall from Section~\ref{sec:intro model} that $X^{(m)}$ is the smallest element for the convex order among $\Z_+$-valued random variables with mean $m$, as proved later in Lemma~\ref{lem:Xm minimal}. By Remark~\ref{rem:N(m)}, it is supported by $\{\lfloor m\rfloor,\lceil m\rceil\}\subset\{0,\dots,N(m)\}$ when $m\geq m_1$.

\begin{corollary}\label{cor:maximal percolation threseshold}
For every $m\geq m_1\approx 1.5682$ and every $\Z_+$-valued random variable $Y$ with mean $m$ whose law differs from that of $X^{(m)}$, we have $p_c(X^{(m)}) > p_c(Y)$.
In particular, for every integer $m\geq 2$, the $m$-ary tree is the unique Bienaym\'e tree with mean $m$ whose grandparent graph has the largest critical parameter.
Moreover, as $m\to\infty$,
\[
\sup\enstq{p_c(Y)}{\Ec{Y}=m}=p_c(X^{(m)})=\frac{1}{m^2}-\frac{1}{m^3}+O\left(\frac{1}{m^4}\right).
\]
\end{corollary}
 The corollary thus shows that the upper bound $p_c(Y)\leq\frac1{m^2}$ of Section~\ref{sec:intro model}, valid for every $Y$ with mean $m$, is sharp at first order as $m\to\infty$.

\subsection{The assumptions of \texorpdfstring{Theorem~\ref{thm:convex order + bound support + large enough mean = comparison}}{Theorem \getrefnumber{thm:convex order + bound support + large enough mean = comparison}} cannot be removed}
The next result ensures that convex ordering alone is not enough to compare the critical parameters associated to two different distributions, which justifies why we added some extra assumption in Theorem~\ref{thm:convex order + bound support + large enough mean = comparison}.
\begin{theorem} \label{thm:contre exemple ordre convexe support}
	Let $m>1$.
	\begin{enumerate}[label=(\roman*)]
	 \item\label{it:counter-example unbounded support} For any $X$ with mean $m$ and unbounded support, there exists $Y$ with
	 \begin{equation*}
	 		 X \prec_{\mathrm{cx}} Y \qquad \text{and} \qquad	p_c(X) < p_c(Y).
	 \end{equation*}
	 \item \label{it:counter-example bounded support}  There exists a random variable $X$ with mean $m$ and support bounded above by an explicit integer $N'(m)$, see \eqref{eq:definition N'}, satisfying $N'(m)\leq (2+o(1))\,m^4$ as $m\to\infty$, and a random variable $Y$ such that
	 \begin{equation*}
	X \prec_{\mathrm{cx}} Y \qquad \text{and} \qquad	p_c(X) < p_c(Y).
\end{equation*}
	\end{enumerate}
\end{theorem}
For instance, we have $N'(2)=84$ by Remark~\ref{rem:N'(2) numerics}, to be compared with $N(2)=5$; the optimal value of $N(m)$ in Theorem~\ref{thm:convex order + bound support + large enough mean = comparison} is therefore between $(1+o(1))m^4$ and $(2+o(1))m^4$. 

Now, recall that the result of Theorem~\ref{thm:convex order + bound support + large enough mean = comparison} and Corollary~\ref{cor:maximal percolation threseshold} only informs us in the case $m\geq m_1 \approx 1.5682$.  
The next result tells us that, even though the threshold $m_1$ is probably not sharp, one needs some assumption of this type. 
\begin{theorem}\label{thm:contre-exemple m proche de 1}
	There exists $\delta>0$ such that for any $m\in\intervalleoo{1}{1+\delta}$, there exists a $\Z_+$-valued random variable $Y$ with mean $m$ such that
	\begin{align*}
		p_c(X^{(m)}) < p_c(Y).
	\end{align*}
\end{theorem}
Since the law of $X^{(m)}$ is minimal for the convex order among distributions with mean $m$, such a $Y$ satisfies $X^{(m)}\prec_{\mathrm{cx}}Y$, so that the conclusions of Theorem~\ref{thm:convex order + bound support + large enough mean = comparison} and Corollary~\ref{cor:maximal percolation threseshold} fail when $m$ is close to $1$.

\subsection{How small can \texorpdfstring{$p_c$}{p\_c} be?}
In the end, one might wonder about the minimal value of $p_c$ across all possible distributions with mean $m$.
We will see that, if the offspring distribution may charge $0$, then $p_c$ can be arbitrarily close to $0$. 
However, the laws that we use for this purpose put a mass close to $1$ on $0$, so that the tree $T$ is infinite with a probability which tends to $0$: the event that the root cluster is infinite has positive probability, but this probability is very small, and the result says more about the degeneracy of the tree than about the geometry of $\GF(T)$. 
It is therefore natural to also consider trees without leaves, for which $T$ is almost surely infinite. 
In that case, the answer depends on whether $m\geq2$ or not.

\begin{theorem}\label{thm:contre exemple pc tend vers 0}
	Let $m>1$.
	\begin{enumerate}[label=(\roman*)]
	\item\label{it:pc to 0 with leaves} There exists a sequence $(X_n)_{n \geq 1}$ of $\Z_+$-valued random variables with mean $m$ such that $p_c(X_n) \rightarrow 0$ as $n \rightarrow \infty$.
	\item\label{it:pc to 0 without leaves} If $m\geq 2$, the random variables $X_n$ can be chosen such that $X_n\geq 1$ almost surely (that is, the trees have no leaves).
	\item\label{it:pc without leaves m<2} If $1<m<2$, then the infimum $\inf\enstq{p_c(X)}{X\geq 1 \text{ a.s.},\ \Ec{X}=m}$ is positive, and it is equivalent to $\frac{2-m}{4}$ as $m\uparrow2$.
	\end{enumerate}
\end{theorem}

\subsection{Related models and corresponding results}

Percolation on trees with long edges has mostly been studied in the oriented setting. 
In multirange percolation on oriented regular trees, short edges link each vertex to its children and long edges link it to its descendants $k$ generations below, with different parameters for short and long edges. 
This model was studied in \cite{dLRV19,dLSV23,CGR25}.
In contrast, our model is non-oriented, which is the source of the main difficulty (a vertex can be connected to the root without its parent being so), and also deals with Bienaymé trees instead of regular trees. 

Let us mention that the first construction of a graph of this kind dates back to Trofimov \cite{Tro85} as an example of a transitive non-unimodular graph.
In his case, the object was obtained by adding to the $d$-regular tree an edge between every vertex and its grandparent \emph{along a fixed end} $\xi$. 
This graph is called the ``grandmother graph'' in \cite{Tim06}; see e.g.\ \cite{PPS06,Tim06,Hut20} for its role in percolation theory.
The descendants of a vertex span a copy of $\GF(T)$ for $X=d-1$ almost surely, that is for $X=X^{(d-1)}$, and it is easy to check that both graphs have the same critical parameter $p_c(X^{(d-1)})$.
For $d=3$, this critical parameter was computed numerically by \v{S}pakulov\'a \cite[Section~5.3]{Spa09} with a general method for transitive tree-like graphs, the grandmother graph being called the ``grandparent tree'' there.
The value $0.158656326\dots$ obtained there agrees with the value of $p_c(X^{(2)})$ given in Remark~\ref{rem:binary tree}.

\subsection{Organisation of the paper}
In Section~\ref{sec:determination}, we formally define the model and introduce Markovian exploration of the root cluster. 
This allows us to deduce that $p_c(X)$ is the unique solution of an explicit equation (Theorem~\ref{prop:the critical percolation parameter is the solution to an equation}), from which we derive in Proposition~\ref{prop:asymptotics pc Xm} the asymptotics of $p_c(X^{(m)})$.
In Section~\ref{sec:monotonicity}, we derive from this equation a criterion to compare critical parameters (Proposition~\ref{prop:criterion for comparison of critical points})  and prove Theorem~\ref{thm:convex order + bound support + large enough mean = comparison} and Corollary~\ref{cor:maximal percolation threseshold}. 
Section~\ref{sec:counterexamples} is devoted to the counter-examples of Theorem~\ref{thm:contre exemple ordre convexe support} and Theorem~\ref{thm:contre-exemple m proche de 1}, and to the proof of Theorem~\ref{thm:contre exemple pc tend vers 0} on small critical parameters.
Finally, the appendix compares annealed and quenched critical parameters.

\subsection*{Acknowledgements}
The authors are grateful to Olivier Couronné for interesting discussions on the directed version of this problem and some of its variations, which led to asking the question investigated in this work.
\subsection*{AI use}
The authors used Claude Opus 5.5 to proofread the paper, in terms of both mathematics and language, 
to write proofs of technical statements, and 
to improve the writing of some proofs.

More precisely, Claude wrote the proofs of Proposition~\ref{prop:clusterGW}, of points~\ref{it:p mapsto fXp is continuous decreasing} and~\ref{it:p,v mapsto Fpv is continuous} of Lemma~\ref{lem:proprietes f_X et F_m}, of Proposition~\ref{prop:asymptotics pc Xm}, of Corollary~\ref{cor:maximal percolation threseshold}, and of Lemmas~\ref{lem:Xm minimal} and~\ref{lem:lower bound without leaves}.
Starting from sketches by the authors, Claude completed the proofs of Theorems~\ref{thm:contre exemple ordre convexe support} and~\ref{thm:contre-exemple m proche de 1}.
It helped improve the authors' proofs of Propositions~\ref{prop:exploration} and~\ref{prop:criterion for comparison of critical points}, of Theorem~\ref{prop:the critical percolation parameter is the solution to an equation}, of point~\ref{it:universal bound Fxp} of Lemma~\ref{lem:proprietes f_X et F_m}, and of Lemmas~\ref{lem:strict phi psi} and~\ref{lem:local spread}.

Two nice catches during the proofreading were: 
It found an error in the authors' calculation of the value $N(m)$ in the proof of Lemma~\ref{lem:bound Nm}.
It realized that the original argument for Theorem~\ref{thm:contre exemple pc tend vers 0} for trees without leaves fails when $m<2$, which led to investigating separately the case $m<2$ and $m\geq 2$, yielding points~\ref{it:pc to 0 without leaves} and~\ref{it:pc without leaves m<2} of this theorem.

Finally, Claude drew Figure~\ref{fig:exploration} and computed the numerical values given in Remarks~\ref{rem:binary tree} and~\ref{rem:N'(2) numerics}.

The authors also tried to use Gemini 3.1 Pro at an earlier stage of the project: it did not significantly help towards the main result, but it suggested the main ingredient of the proof of Lemma~\ref{lem:strict phi psi}, and provided a trick to get a better constant in the inequality proven in Lemma~\ref{lem:inequality difference derivatives}.

The authors take full responsibility for the mathematical correctness of the paper.

\section{The model and its critical parameter}\label{sec:determination}

In this section, $X$ is a random variable with law $\mu$ on $\Z_+=\{0,1,2,\dots\}$, expectation $m > 1$ and probability generating function $\phi$.

\subsection{Notation}\label{sec:notation}
We consider here the Ulam--Harris tree, with set of vertices given by
\begin{align*}
	\mathbb U=\bigcup_{n=0}^\infty \{1,2,\dots\}^n,
\end{align*}
where, we recall, $\{1,2,\dots\}^0=\{\emptyset\}$.
For any vertex $u$ we denote by $\haut(u)$ its height in the tree and $\hat{u}$ its parent.
For reasons that will become clear later on, we add a new vertex $\hat{\emptyset}$ to the tree, as a parent of the root $\emptyset$, with height $\haut(\hat\emptyset):=-1$.
We write $u\preceq v$ to denote that $u$ is an ancestor of $v$.
For every $u\in \mathbb{U}$, we define the shift $\sigma_u$, which maps $\hat u$ to $\hat\emptyset$ and $uv$ to $v$ for every $v\in\mathbb U$.

We consider the Bienaym\'e tree $T$ with offspring distribution $\mu$ as a subset of $\mathbb U$ in the usual way, and for $u\in T$ we denote by $T_u:=\enstq{v\in T}{v\succeq u}$ the subtree of the descendants of $u$, including $u$ itself. 
Note that $T$ is finite with positive probability when $\mu(0)>0$.

Now, introduce $\GF(T)$, the associated grandparent graph, as the graph with vertex set $T \cup \{\hat{\emptyset}\}$ and set of edges given by
\begin{align*}
\enstq{(\hat{u},u)}{u \in T} \cup \enstq{(\doublehat{u},u)}{u \in T \setminus \{\emptyset\}},
\end{align*}
which in other words is just saying that every vertex is connected by an edge to its parent, and to its grandparent.
We will sometimes refer to the first type of edges (parent--child) as \emph{short} edges, and to the second type of edges (grandparent--child) as \emph{long} edges.
Conditionally on $T$, we assign a random variable $U_e$ uniform on $\intervalleff{0}{1}$ on every edge $e$ of $\GF(T)$, independently for every edge $e$.

Fix $p\in\intervalleff{0}{1}$. 
We perform unoriented bond percolation on $\GF(T)$ by retaining every edge $e$ of $\GF(T)$ if and only if $U_e<p$. 
We call the retained edges open and the remaining edges closed. 
We denote by $\mathcal{C}_p$ the set of vertices of $T$ which are connected to $\hat{\emptyset}$ by a path of open edges, that is the open cluster of $\hat\emptyset$ without the vertex $\hat\emptyset$ itself. 
It is clear from the monotonicity in $p$ for inclusion of the above construction that the mapping $p\mapsto \Pp{|\mathcal{C}_p| = \infty}$ is non-decreasing.

Note that the definition of $\GF(T)$ here differs from the one given in the introduction through the addition of the extra vertex $\hat{\emptyset}$.
As we will see in the next subsection, this will make some cluster exploration arguments easier to write.
This addition does not change the critical parameter: the edges incident to $\hat\emptyset$ do not belong to the graph of the introduction, and $\hat\emptyset$ has finitely many neighbours, so that $|\mathcal C_p|=\infty$ with positive probability if and only if the open cluster of $\emptyset$ in the graph of the introduction is infinite with positive probability.
Hence the critical parameter defined in the introduction is
\begin{align*}
	p_c(X) = \inf\enstq{p \in \intervalleff{0}{1}}{\Pp{|\mathcal{C}_p| = \infty} >0},
\end{align*}
so that $\Pp{|\mathcal{C}_p| = \infty}>0$ if $p>p_c(X)$ and $\Pp{|\mathcal{C}_p| = \infty}=0$ if $p<p_c(X)$.

\subsection{Exploring the root cluster}
\begin{figure}
\centering
\begin{tikzpicture}[scale=0.8]
 \gfcoords
 \foreach \h in {-1,...,3} {\draw[gray!25] (-3.8,\h) -- (2.8,\h); \node[left,font=\scriptsize] at (-3.8,\h) {$\h$};}
 \node[font=\scriptsize,above] at (-4.05,3.25) {height};
 \foreach \u in {a2,b1,a12,c} {\subtreeat{\u}}
 \foreach \x/\y in {r0/r,r/a,r/b,r/c,a/a1,a/a2,b/b1,a1/a11,a1/a12} {\draw[shortedge] (\x) -- (\y);}
 \draw[openlong] (r0) to[bend left=30] (a);
 \draw[openlong] (r0) to[bend right=35] (b);
 \draw[closedlong] (r0) to[bend right=30] (c);
 \draw[openlong] (r) to[bend right=15] (a1);
 \draw[closedlong] (r) to (a2);
 \draw[closedlong] (r) to[bend left=35] (b1);
 \draw[openlong] (a) to[bend left=25] (a11);
 \draw[closedlong] (a) to (a12);
 \foreach \u in {r,a,b,a1,a11} {\node[vin] at (\u) {};}
 \foreach \u in {c,a2,b1,a12} {\node[vbd] at (\u) {};}
 \node[vroot] at (r0) {};
 \node[below,font=\small] at (r0) {$\hat\emptyset$};
 \node[right,font=\small] at (r) {$\emptyset$};
\end{tikzpicture}\hfill
\begin{tikzpicture}[scale=0.8]
 \gfcoords
 \draw[openlong] (r0) to[bend left=30] (a);
 \draw[openlong] (r0) to[bend right=35] (b);
 \draw[openlong] (r) to[bend right=15] (a1);
 \draw[openlong] (a) to[bend left=25] (a11);
 \foreach \x/\y in {r0/r,r/a,r/b,a/a1,a1/a11} {\draw[shortdash] (\x) -- (\y);}
 \foreach \u in {a,b,a11} {\node[vin] at (\u) {};}
 \foreach \u in {r,a1} {\node[vev] at (\u) {};}
 \node[vroot] at (r0) {};
 \node[below,font=\small] at (r0) {$\hat\emptyset$};
 \node[right,font=\small] at (r) {$\emptyset$};
 \node at (2.1,1) {}; \node at (-3.8,3) {};
\end{tikzpicture}\hfill
\begin{tikzpicture}[scale=0.8]
 \gfcoords
 \foreach \u in {a2,b1} {\subtreeat{\u}}
 \draw[openlong] (r0) to[bend left=30] (a);
 \draw[openlong] (r0) to[bend right=35] (b);
 \draw[openlong] (a) to[bend left=25] (a11);
 \draw[shortedge] (a) -- (a2); \draw[shortedge] (b) -- (b1);
 \foreach \u in {r,a1} {\node[vghost] at (\u) {};}
 \foreach \u in {a,b,a11} {\node[vin] at (\u) {};}
 \foreach \u in {a2,b1} {\node[vbd] at (\u) {};}
 \node[vroot] at (r0) {};
 \node[below,font=\small] at (r0) {$\hat\emptyset$};
 \node[right,font=\small,gray] at (r) {$\emptyset$};
 \node[right,font=\small] at (a2) {$v$}; \node[left,font=\small] at (a) {$\hat v$};
 \node[right,font=\small] at (b1) {$v'$}; \node[right,font=\small] at (b) {$\hat v'$};
 \node[font=\small] at (-1.1,3.2) {$T_v$}; \node[font=\small] at (0.7,3.2) {$T_{v'}$};
 \node at (2.1,1) {}; \node at (-3.8,3) {};
\end{tikzpicture}
\caption{Exploration of the root cluster. Left: a portion of $\GF(T)$, with short edges in grey. The vertices of $\widetilde{\mathcal C_p}$, drawn as black discs, are linked to their grandparents by open long edges, drawn as solid lines; the vertices of $\partial\widetilde{\mathcal C_p}$, drawn as white discs, are linked to theirs by closed long edges, drawn as dotted lines; the grey triangles represent the subtrees $T_u$ for $u\in\partial\widetilde{\mathcal C_p}$. Middle: the set $\widetilde{\mathcal C_p}\cup\{\hat\emptyset\}$, with $\widetilde{\mathcal C_p}^o$ in black and $\widetilde{\mathcal C_p}^e$ in grey. Open long edges preserve the parity of the height, and the only other edges between these vertices are the $|\widetilde{\mathcal C_p}|=5$ short edges $(\hat u,u)$ with $u\in\widetilde{\mathcal C_p}$, drawn as dashed lines, which link the two parity classes. Right: if these five short edges are closed, then $\mathcal D_p=\widetilde{\mathcal C_p}^o$ and $\partial\mathcal D_p=\{v,v'\}$. The subtree $T_v$ is attached to the rest of the open graph only through $\hat v$, so that the exploration of $\mathcal C_p$ continues independently in $\{\hat v\}\cup T_v$ and in $\{\hat v'\}\cup T_{v'}$.}
\label{fig:exploration}
\end{figure}
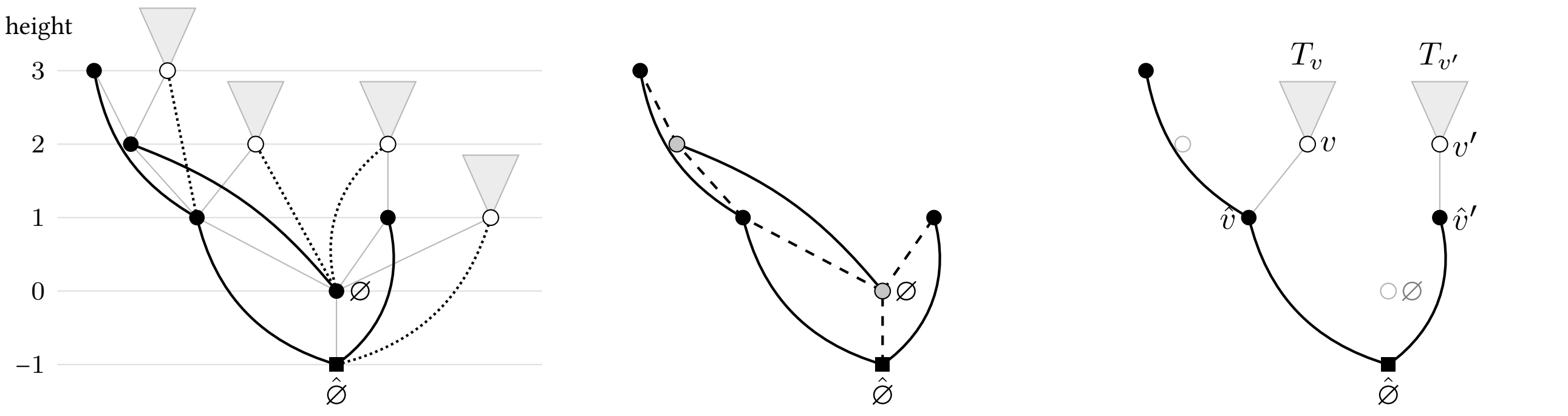

The difficulty of studying a non-oriented version of percolation is that it can be the case that a vertex is connected to the root without its parent being in the cluster as well.
This makes exploration arguments harder to deal with.

We introduce the following definition: for a vertex $u\in T\setminus\{\emptyset\}$ and its parent $\hat{u}$, we write $\hat{u} \overset{p}{\rightsquigarrow} u$ if the long edge $(\doublehat{u},u)$ is open. 
For any vertex $u$, we also set $u\overset{p}{\rightsquigarrow} u$ by definition.
We extend this relation by taking its transitive closure, which we still denote by $\overset{p}{\rightsquigarrow}$.
In other words, $\emptyset\overset{p}{\rightsquigarrow}u$ if and only if the long edges $(\doublehat{w},w)$ are open for every vertex $w\neq\emptyset$ on the path from $\emptyset$ to $u$ in $T$.

We now introduce two variants of the cluster, together with their boundaries:
\begin{align*}
	&\widetilde{\mathcal{C}_p} := \enstq{u \in T}{\emptyset \overset{p}{\rightsquigarrow} u} \qquad &\text{and} \qquad &\partial \widetilde{\mathcal{C}_p} := \enstq{u \in T}{\emptyset \overset{p}{\rightsquigarrow} \hat{u} \text{ and } \ \emptyset \overset{p}{\nrightsquigarrow} u},\\
	&\mathcal{D}_p := \mathcal{C}_p \cap \widetilde{\mathcal{C}_p} \qquad &\text{and} \qquad &\partial \mathcal{D}_p := \enstq{u \in \partial \widetilde{\mathcal{C}_p} }{\hat{u} \in \mathcal{D}_p }.
\end{align*}
In words, $\partial\widetilde{\mathcal C_p}$ is the set of children of vertices of $\widetilde{\mathcal C_p}$ whose long edge to their grandparent is closed; in particular $\hat u\in\widetilde{\mathcal C_p}$ for every $u\in\partial\widetilde{\mathcal C_p}$. 
We denote by $\partial^e\widetilde{\mathcal C_p}$ and $\partial^o\widetilde{\mathcal C_p}$ the sets of vertices of $\partial\widetilde{\mathcal C_p}$ with even and odd height respectively, and by $\widetilde{\mathcal{C}_p}^e$ and $\widetilde{\mathcal{C}_p}^o$ the sets of vertices of $\widetilde{\mathcal{C}_p}$ with even and odd height respectively.
Finally, for $u\in T$, we consider the set of vertices of $T_u$ which are connected to $\hat u$ in the subgraph of $\GF(T)$ induced by $\{\hat u\}\cup T_u$, namely
\[
\mathcal C_p(u):=\enstq{v\in T_u}{\hat u \text{ and } v \text{ are connected by an open path of } \GF(T) \text{ staying in } \{\hat u\}\cup T_u}.
\]
Note that $\mathcal C_p(\emptyset)=\mathcal C_p$.
The main result of this section is a Markovian exploration of the root cluster.
\begin{proposition}\label{prop:exploration}
Let $p\in\intervalleoo{0}{\frac1m}$. 
The sets $\widetilde{\mathcal C_p}$, $\partial\widetilde{\mathcal C_p}$, $\mathcal{D}_p$ and $\partial \mathcal{D}_p$ are finite almost surely, and
\[
	\mathcal{C}_p = \mathcal{D}_p \sqcup \bigsqcup_{v \in \partial \mathcal{D}_p} \mathcal{C}_p(v).
\]
Moreover, conditionally on $\mathcal{D}_p $ and $\partial \mathcal{D}_p $, the family
$\left(\sigma_v(\mathcal{C}_p(v)) \right)_{v \in \partial \mathcal{D}_p }$ is i.i.d.\ with the same law as $\mathcal{C}_p  $.
\end{proposition}

\begin{proof}
First, notice that $\widetilde{\mathcal{C}_p}$ is a Bienaym\'e tree whose offspring distribution is the law of $\sum_{i=1}^{X}\varepsilon_i$, where $(\varepsilon_i)_{i\geq1}$ are i.i.d.\ Bernoulli random variables with parameter $p$, independent of $X$: indeed, a child of a vertex of $\widetilde{\mathcal C_p}$ belongs to $\widetilde{\mathcal C_p}$ if and only if the long edge linking it to its grandparent is open. 
This offspring distribution has mean $pm$. 
Hence, since $pm<1$, it is almost surely finite, with total size of finite expectation. 
Since every vertex of $\widetilde{\mathcal C_p}$ has on average $(1-p)m$ children in $\partial\widetilde{\mathcal C_p}$, we also have $\Ec{|\partial\widetilde{\mathcal C_p}|}=(1-p)m\,\Ec{|\widetilde{\mathcal C_p}|}<\infty$. 
Since $\mathcal D_p\subset\widetilde{\mathcal C_p}$ and $\partial\mathcal D_p\subset\partial\widetilde{\mathcal C_p}$, all four sets are finite almost surely.

\step{description of $\mathcal D_p$}\label{step:description Dp} This step is illustrated in the middle of Figure~\ref{fig:exploration}.
Every $u\in\widetilde{\mathcal C_p}\setminus\{\emptyset\}$ is linked by an open long edge to its grandparent $\doublehat u\in\widetilde{\mathcal C_p}\cup\{\hat\emptyset\}$, which has the same parity. Following these edges downwards, we see that every vertex of $\widetilde{\mathcal{C}_p}^o$ is connected to $\hat\emptyset$ and every vertex of $\widetilde{\mathcal{C}_p}^e$ is connected to $\emptyset$. Hence $\widetilde{\mathcal{C}_p}^o\subset\mathcal{C}_p$, and either $\widetilde{\mathcal{C}_p}^e\subset\mathcal{C}_p$ or the sets $\widetilde{\mathcal{C}_p}^e$ and $\mathcal{C}_p$ are disjoint.
Recall that $\widetilde{\mathcal C_p}$ spans a subtree of $T$ containing $\emptyset$. Hence, apart from these long edges, the only edges of $\GF(T)$ between vertices of $\widetilde{\mathcal C_p}\cup\{\hat\emptyset\}$ are the $|\widetilde{\mathcal{C}_p}|$ short edges $(\hat u,u)$ for $u\in\widetilde{\mathcal C_p}$, and each of them links a vertex of $\widetilde{\mathcal{C}_p}^e$ to a vertex of $\widetilde{\mathcal{C}_p}^o\cup\{\hat\emptyset\}$. 
If one of them is open, then $\widetilde{\mathcal{C}_p}^e\subset\mathcal{C}_p$.

Conversely, the two parity classes cannot be connected through the rest of the graph, because of the following observation (see the right of Figure~\ref{fig:exploration}):
\begin{equation}\label{eq:cut vertex}
\text{for } u\in\partial\widetilde{\mathcal C_p},\ \text{every open edge between } T_u \text{ and its complement has } \hat u \text{ as an endpoint.}
\end{equation}
Indeed, the edges of $\GF(T)$ between $T_u$ and its complement are $(\hat u,u)$, $(\doublehat u,u)$ and the long edges $(\hat u,w)$ for $w$ a child of $u$, and $(\doublehat u,u)$ is closed by definition of $\partial\widetilde{\mathcal C_p}$. Moreover, every vertex of $T\setminus\widetilde{\mathcal C_p}$ has a unique ancestor in $\partial\widetilde{\mathcal C_p}$, namely its first ancestor which is not in $\widetilde{\mathcal C_p}$, so that
\begin{equation}\label{eq:partition of T}
T=\widetilde{\mathcal C_p}\sqcup\bigsqcup_{u\in\partial\widetilde{\mathcal C_p}}T_u .
\end{equation}

By~\eqref{eq:cut vertex}, an open path between vertices in $\widetilde{\mathcal C_p}\cup\{\hat\emptyset\}$ which enters some $T_u$ with $u\in\partial\widetilde{\mathcal C_p}$ can only leave it through the vertex $\hat u$ by which it entered. 
Removing all such excursions ``above'' $\widetilde{\mathcal C_p}\cup\{\hat\emptyset\}$ yields an open path that stays within $\widetilde{\mathcal C_p}\cup\{\hat\emptyset\}$. 
Hence two vertices of $\widetilde{\mathcal C_p}\cup\{\hat\emptyset\}$ that are connected by an open path are connected by an open path using only vertices of $\widetilde{\mathcal C_p}\cup\{\hat\emptyset\}$. 
	In particular, if all the short edges $(\hat u,u)$ with $u\in\widetilde{\mathcal C_p}$ are closed, then the sets $\widetilde{\mathcal{C}_p}^e$ and $\mathcal{C}_p$ are disjoint.
To sum up:
\begin{itemize}
	\item either one of these short edges is open, and $\mathcal{D}_p=\widetilde{\mathcal C_p}$,
	\item or all of them are closed, and $\mathcal{D}_p=\widetilde{\mathcal{C}_p}^o$.
\end{itemize}

\step{decomposition}\label{step:decomposition} Let $u\in\partial\widetilde{\mathcal C_p}$. By~\eqref{eq:cut vertex}, the set $\mathcal C_p\cap T_u$ is non-empty only if $\hat u\in\mathcal C_p$; since $\hat u\in\widetilde{\mathcal C_p}$, this happens if and only if $\hat u\in\mathcal D_p$, that is $u\in\partial\mathcal D_p$. In that case, an open path from $\hat\emptyset$ to a vertex of $T_u$ stays in $T_u$ after its last visit to $\hat u$, again by~\eqref{eq:cut vertex}, so that $\mathcal C_p\cap T_u=\mathcal C_p(u)$. 
Since $\mathcal C_p\subset T$ and $\mathcal C_p\cap\widetilde{\mathcal C_p}=\mathcal D_p$, intersecting the partition~\eqref{eq:partition of T} with $\mathcal C_p$ gives the decomposition of $\mathcal C_p$.

\step{independence}\label{step:independence} This is a consequence of the branching property. Let $S$ and $S'$ be finite subsets of $\mathbb U$. The event $\{\widetilde{\mathcal C_p}=S,\ \partial\widetilde{\mathcal C_p}=S'\}$ and the states of the short edges $(\hat u,u)$ for $u\in S$, which determine $\mathcal D_p$ and $\partial\mathcal D_p$ by Step~\ref{step:description Dp}, only depend on the numbers of children of the vertices of $S$, on the states of these short edges, and on the states of the long edges $(\doublehat u,u)$ for $u\in(S\cup S')\setminus\{\emptyset\}$. On the other hand, for $v\in S'$, the set $\mathcal C_p(v)$ only depends on the subtree $T_v$ and on the states of the edges of $\GF(T)$ with both endpoints in $\{\hat v\}\cup T_v$. None of these edges and none of the vertices of $T_v$ are involved in the event above, and these sets of edges are disjoint for distinct $v\in S'$. By the branching property of Bienaym\'e trees, conditionally on this event, the subtrees $(T_v)_{v\in S'}$ together with the states of these edges are i.i.d.\ and, after applying the shift $\sigma_v$, have the same law as $T$ together with the states of the edges of $\GF(T)$. Hence, conditionally on $(\mathcal D_p,\partial\mathcal D_p)$, the sets $\sigma_v(\mathcal C_p(v))$ for $v\in\partial\mathcal D_p$ are i.i.d.\ with the law of $\mathcal C_p$.
\end{proof}

As a byproduct of the previous proof we have the following criterion for percolation.

\begin{proposition}\label{prop:clusterGW}
Let $p\in\intervalleoo{0}{\frac1m}$. We have
\begin{align*}
	\Pp{|\mathcal{C}_p | = \infty} >0 \quad \iff \quad \Ec{|\partial \mathcal{D}_p  |}>1,
\end{align*}
and
\[
\Ec{|\partial \mathcal{D}_p  |} = \Ec{|\partial^e \widetilde{\mathcal{C}_p}  |} +
\Ec{\left(1-(1-p)^{|\widetilde{\mathcal{C}_p} |} \right)|\partial^o \widetilde{\mathcal{C}_p}  |}.
\]
\end{proposition}
\begin{proof}
Iterating Proposition~\ref{prop:exploration}, the cluster $\mathcal C_p$ is the disjoint union of finite ``blocks'' (the first block being $\mathcal D_p$), which are organised as a Bienaym\'e tree whose offspring distribution is the law of $|\partial\mathcal D_p|$.
A block may be empty, but then it has no child, since $\partial\mathcal D_p$ is empty when $\mathcal D_p$ is.
Hence $|\mathcal C_p|=\infty$ if and only if this Bienaym\'e tree is infinite, which has positive probability if and only if $\Ec{|\partial\mathcal D_p|}>1$ or $|\partial\mathcal D_p|=1$ almost surely. 
The latter case is excluded: with positive probability, all the children of $\emptyset$ have a closed long edge to $\hat\emptyset$, the short edge $(\hat\emptyset,\emptyset)$ is open and $\emptyset$ has at least two children (recall that $m>1$); on this event $\mathcal D_p=\widetilde{\mathcal C_p}=\{\emptyset\}$ and $|\partial\mathcal D_p|\geq 2$.

For the formula, let $u\in\partial\widetilde{\mathcal C_p}$, so that $\hat u\in\widetilde{\mathcal C_p}$. If $u$ has even height, then $\hat u\in\widetilde{\mathcal C_p}^o\subset\mathcal C_p$, so that $u\in\partial\mathcal D_p$. If $u$ has odd height, then by Step~\ref{step:description Dp} of the proof of Proposition~\ref{prop:exploration}, we have $\hat u\in\mathcal D_p$ if and only if one of the $|\widetilde{\mathcal C_p}|$ short edges linking a vertex of $\widetilde{\mathcal C_p}$ to its parent is open. 
Conditionally on $(\widetilde{\mathcal C_p},\partial\widetilde{\mathcal C_p})$, these short edges are independent and open with probability $p$, since their states are not involved in the definition of $\widetilde{\mathcal C_p}$ and $\partial\widetilde{\mathcal C_p}$. 
The formula follows.
\end{proof}
\subsection{A study of the generating function \texorpdfstring{$A$}{A}}

Motivated by Proposition~\ref{prop:clusterGW},
for $p< \frac{1}{m}$ fixed and for any $x,y,z \in \intervalleff{0}{1}$ we introduce
 \begin{align*}
 	A(x,y,z,p):= \Ec{x^{|\widetilde{\mathcal{C}_p}|} \  y^{|\partial^e\widetilde{\mathcal{C}_p}|} \  z^{|\partial^o\widetilde{\mathcal{C}_p}|}}.
 \end{align*}

\begin{proposition}\label{prop:generating function A}
Recall that $\phi$ is the probability generating function of $X$.
The function
\begin{align*}
	A : \quad \intervalleff{0}{1}^3\times \intervallefo{0}{1/m} &\to \mathbb R\\
	(x,y,z,p) &\mapsto A(x,y,z,p)
\end{align*}
is well defined and satisfies the equation
\begin{align}\label{eq:equation satisfied by A}
	A(x,y,z,p)&=x \cdot \phi\left( p A(x,z,y,p) + (1-p) z\right).
\end{align}
For $p<\frac{1}{m}$, we have
\begin{align}\label{eq:taille moyenne du cluster de l'origine}
\Ec{|\partial \mathcal{D}_p  |}=\frac{(1-p)m}{1-pm} - \partial_zA(1-p,1,1,p).
\end{align}
\end{proposition}
\begin{proof}
We already mentioned that, for $p < 1/m$, the random variables appearing in the definition of $A$ are finite almost surely. 
As a consequence, the function $A$ is well defined on $\intervalleff{0}{1}^3\times \intervallefo{0}{1/m}$. 
Decomposing $T$ according to the children of $\emptyset$ we have
\begin{equation}
\begin{pmatrix}
|\widetilde{\mathcal C_p}|\\
|\partial^e \widetilde{\mathcal C_p}|\\
|\partial^o \widetilde{\mathcal C_p}|
\end{pmatrix}
=
\begin{pmatrix}
1\\
0\\
0
\end{pmatrix}
+
\sum_{v \text{ child of } \emptyset}
\begin{pmatrix}
|\widetilde{\mathcal C_p} \cap T_v|\\
|\partial^e \widetilde{\mathcal C_p} \cap T_v|\\
|\partial^o \widetilde{\mathcal C_p} \cap T_v|
\end{pmatrix}.
\end{equation}
Take $v$ a child of $\emptyset$ in $T$. 
If the long edge connecting it to $\hat \emptyset$ is closed, then $v$ belongs to $\partial^o \widetilde{\mathcal C_p}$ and its contribution in the sum is $(0,0,1)^\top$. 
This happens with probability $1-p$.
Conversely, if the long edge connecting it to $\hat \emptyset$ is open, then $v$ belongs to $\widetilde{\mathcal C_p}$.
Conditionally on this long edge being open, the part of $\widetilde{\mathcal C_p}$ above $v$, more precisely $\sigma_v( \widetilde{\mathcal C_p} \cap T_v)$, has the same distribution as $ \widetilde{\mathcal C_p}$.
Note that the parity of generations is swapped in $\sigma_v( \widetilde{\mathcal C_p} \cap T_v)$ compared to $\widetilde{\mathcal C_p}$. 
This gives, for every child $v$ of $\emptyset$,
\[
	\Ecsq{x^{|\widetilde{\mathcal{C}_p} \cap T_v|} \  y^{|\partial^e\widetilde{\mathcal{C}_p}\cap T_v|} \  z^{|\partial^o\widetilde{\mathcal{C}_p}\cap T_v|}}{\text{number of children of } \emptyset}
	=
	(1-p) z + p A(x,z,y,p).
\]
The independence between disjoint subtrees, and the fact that the number of children of $\emptyset$ has the law of $X$, gives
\[
	A(x,y,z,p) = x \mathbb E \left[ \left((1-p) z + p A(x,z,y,p) \right)^{X} \right] = x \phi \left((1-p) z + p A(x,z,y,p) \right).
\]

We turn to the proof of~\eqref{eq:taille moyenne du cluster de l'origine}. Proposition~\ref{prop:clusterGW} directly gives
\[
	\mathbb E \left[ |\partial \mathcal D_p| \right] = \partial_yA(1,1,1,p) + \partial_z A(1,1,1,p) - \partial_z A(1-p,1,1,p),
\]
where the partial derivatives in $y$ and $z$ at $y=1$ or $z=1$ are left derivatives; they are finite and can be computed by differentiating under the expectation (by monotone convergence), since $\Ec{|\partial\widetilde{\mathcal C_p}|}<\infty$.
Differentiating~\eqref{eq:equation satisfied by A} we get
\begin{align*}
	\partial_yA(x,y,z,p)  &= x p\partial_z A(x,z,y,p) \phi'\left( p A(x,z,y,p) + (1-p) z\right),\\
	\partial_z A(x,y,z,p) &= x \left(p \partial_y A(x,z,y,p) +(1-p) \right) \phi'\left( p A(x,z,y,p) + (1-p) z\right).
\end{align*}
Recall that $A(1,1,1,p)=1$. 
Evaluating these two identities at the point $(1,1,1,p)$ and summing, we get
\[
	\partial_yA(1,1,1,p) + \partial_zA(1,1,1,p) = \left( 1-p + p (\partial_yA(1,1,1,p) + \partial_zA(1,1,1,p)) \right) \, \phi'(1).
\]
The result follows with $\phi'(1)=m$.
\end{proof}

The theorem below will pinpoint the critical parameter $p_c(X)$ as the solution to an equation.
We first introduce some notation.
For $p\in\intervallefo{0}{1}$, we define $f_X(p)$ as the unique solution in $\intervalleff{0}{1}$ of the fixed point equation
\begin{align}\label{eq:equation f_X(p)}
	x = 1 - p + p(1-p) \phi(x).
\end{align}
Indeed, for $p=0$ this solution is $1$. 
For $p\in\intervalleoo{0}{1}$, the function $x\mapsto 1-p+p(1-p)\phi(x)-x$ is convex, its value at $0$ is $1-p+p(1-p)\mu(0)>0$ and its value at $1$ is $-p^2<0$, so that it has a unique zero in $\intervalleff{0}{1}$.
Since the right-hand side of~\eqref{eq:equation f_X(p)} takes its values in $\intervalleff{1-p}{1-p^2}$, we have $f_X(p)\in\intervalleff{1-p}{1-p^2}$.
We also define
\begin{equation} \label{eq:definitionF_m}
	F_m(p,v) := (1-p)^2 \cdot \left(1-mp +  p^2 \cdot (m-1) \cdot  v \right) \cdot v - (m-1).
\end{equation}
The next theorem characterises $p_c(X)$ in terms of $f_X$ and $F_m$.
\begin{theorem}\label{prop:the critical percolation parameter is the solution to an equation}
The function $p\mapsto F_m\left(p,\phi'(f_X(p))\right)$ is continuous and strictly decreasing on the interval $\intervalleff{0}{\frac{1}{m^2}}$, with positive value at $0$ and negative value at $\frac{1}{m^2}$.
The critical parameter $p_c(X)$ is the unique solution of
\begin{align} \label{eq:caracterisep_c(X)}
F_m(p,\phi'(f_X(p))) = 0,
\end{align}
on the interval $\intervalleoo{0}{\frac{1}{m^2}}$.
\end{theorem}
We split the proof of this result into two lemmas.
\begin{lemma}\label{lem:p is supercritical iff F<0}
	For any $p\in \intervallefo{0}{\frac{1}{m}}$ we have
\begin{align*}
\Pp{|\mathcal{C}_p | = \infty} >0 \quad \iff \quad F_m(p,\phi'\left(f_X(p) \right)) <0.
\end{align*}
\end{lemma}
\begin{proof}
For $p=0$, both sides of the equivalence are false, since the set $\mathcal C_0$ is empty and $F_m(0,\phi'(f_X(0)))=F_m(0,m)=1$. 
Let now $p\in\intervalleoo{0}{\frac1m}$.
Thanks to Proposition~\ref{prop:clusterGW} and \eqref{eq:taille moyenne du cluster de l'origine} we know that
\begin{align*}
\Pp{|\mathcal{C}_p | = \infty} >0 \quad \underset{\text{Prop.~\ref{prop:clusterGW}}}{\iff} \quad \Ec{|\partial \mathcal{D}_p  |}>1  \underset{\eqref{eq:taille moyenne du cluster de l'origine}}{\iff}  \partial_zA(1-p,1,1,p) < \frac{m-1}{1-pm}.
\end{align*}
We start by studying the quantity $\partial_zA(1-p,1,1,p)$ that appears in the necessary and sufficient condition above.
First, we consider the equation \eqref{eq:equation satisfied by A} satisfied by  $A$ and we
differentiate in $z$ and take $x=1-p$ and $y=z=1$.
We get
\begin{align*}
	\partial_z A(1-p,1,1,p)= (1-p) \cdot \left( p\partial_yA(1-p,1,1,p) + (1-p)\right) \cdot \phi'\left( p A(1-p,1,1,p) + (1-p)\right).
\end{align*}
Similarly, differentiating in $y$ and taking $x=1-p$ and $y=z=1$ in \eqref{eq:equation satisfied by A}, we get
\begin{align*}
		\partial_y A(1-p,1,1,p)&=(1-p)\cdot\left( p\partial_zA(1-p,1,1,p)\right) \cdot  \phi'\left( p A(1-p,1,1,p) + (1-p)\right).
\end{align*}
Plugging the expression for $\partial_y A(1-p,1,1,p)$ given in the last display into the right-hand side of the previous one yields
\begin{multline} \label{eq:calcul d_z A}
	\partial_z A(1-p,1,1,p)
	= (1-p) \cdot \phi'\left( p A(1-p,1,1,p) + (1-p)\right) \\
	\cdot \left( p (1-p)\cdot\left( p\partial_zA(1-p,1,1,p)\right) \cdot  \phi'\left( p A(1-p,1,1,p) + (1-p)\right) + (1-p)\right).
\end{multline}
Using~\eqref{eq:equation satisfied by A} with $x=1-p$ and $y=z=1$, we see that
\begin{align*}
A(1-p,1,1,p) = (1-p) \phi\left( p A(1-p,1,1,p) + (1-p)\right)
\end{align*}
so that the quantity $	p A(1-p,1,1,p) + (1-p)$ satisfies
\begin{align*}
	p A(1-p,1,1,p) + (1-p) = 1-p + p(1-p) \phi\left( p A(1-p,1,1,p) + (1-p)\right)
\end{align*}
which corresponds to the fixed point equation \eqref{eq:equation f_X(p)}, so that $p A(1-p,1,1,p) + (1-p) = f_X(p)$ by uniqueness of the fixed point.

Now, starting from \eqref{eq:calcul d_z A}, we substitute
 $p A(1-p,1,1,p) + (1-p)$ by $f_X(p)$.
 This yields
\begin{align*}
\partial_z A(1-p,1,1,p)
= (1-p) \cdot \phi'\left( f_X(p) \right)
\cdot \left( p (1-p)\cdot\left( p \partial_z A(1-p,1,1,p)  \right) \cdot  \phi'\left(f_X(p) \right) + (1-p)\right),
\end{align*}
which can be re-written as
\begin{align*}
	\partial_z A(1-p,1,1,p)
	= \frac{(1-p)^2\cdot \phi'\left( f_X(p) \right) }{1 - (1-p)^2p^2 (\phi'\left(f_X(p) \right))^2 },
\end{align*}
where the denominator on the right-hand side is positive because $\phi'\left(f_X(p) \right)\leq \phi'(1) = m$ and $p<\frac{1}{m}$.
Now, $\Pp{|\mathcal{C}_p | = \infty} >0 $ if and only if $\partial_zA(1-p,1,1,p) < \frac{m-1}{1-pm}$, which, using the equality above, is equivalent to
\begin{align*}
	\frac{(1-p)^2\cdot \phi'\left( f_X(p) \right) }{1 - (1-p)^2p^2 (\phi'\left(f_X(p) \right))^2 } < \frac{m-1}{1-pm}.
\end{align*}
After some simplifications of the last display (multiply both sides by the positive quantity $(1-pm)\left(1 - (1-p)^2p^2 \phi'(f_X(p))^2\right)$), we get that
\begin{align*}
\Pp{|\mathcal{C}_p | = \infty} >0 \quad \iff \quad F_m(p,\phi'\left(f_X(p) \right)) <0,
\end{align*}
which is the result stated above.
\end{proof}
\begin{lemma} \label{lem:proprietes f_X et F_m}
We have the following properties on the functions $f_X$ and $F_m$:
\begin{enumerate}[label=(\alph*)]
	\item \label{it:p mapsto fXp is continuous decreasing} The function $p\mapsto f_X(p)$ is continuous and strictly decreasing on $\intervalleff{0}{\frac1m}$.
	\item\label{it:universal bound Fxp} For $p\in\intervalleof{0}{\frac1m}$, we have the following bounds on $f_X(p)$,
	\begin{align*}
		1-p \leq 1 - \frac{p^2}{1-p(1-p)m} \leq f_X(p) < 1 - p^2 \quad
		\text{or equivalently} \quad p^2 < 1-f_X(p) \leq \frac{p^2}{1-p(1-p)m} \leq p.
	\end{align*}
	\item \label{it:p,v mapsto Fpv is continuous} The function $(p,v) \mapsto F_m(p,v)$ is continuous on $\intervallefo{0}{\frac{1}{m}}\times \intervallefo{0}{+\infty}$, strictly increasing in $v$ for fixed $p\in \intervallefo{0}{\frac{1}{m}}$ and strictly decreasing in $p\in\intervalleff{0}{\frac1{m^2}}$ for fixed $v\in \intervalleof{0}{m}$.
\end{enumerate}
\end{lemma}
\begin{proof}

We first prove \ref{it:p mapsto fXp is continuous decreasing}. 
Let $H(p,x):=1-p+p(1-p)\phi(x)-x$, which is $C^1$ on $\R\times\intervalleoo{0}{1}$.
For $p\in\intervalleof{0}{\frac1m}$, we have $f_X(p)\in\intervalleff{1-p}{1-p^2}\subset\intervalleoo{0}{1}$ and, since $\phi'\leq m$ on $\intervalleff{0}{1}$ and $pm\leq1$,
\[
\partial_xH(p,f_X(p))=p(1-p)\phi'(f_X(p))-1\leq p(1-p)m-1\leq -p<0.
\]
By the implicit function theorem, $f_X$ is $C^1$ on $\intervalleof{0}{\frac1m}$ with
\[
f_X'(p)=-\frac{\partial_pH(p,f_X(p))}{\partial_xH(p,f_X(p))},\qquad \partial_pH(p,x)=-1+(1-2p)\phi(x)\leq -1+|1-2p|<0 .
\]
Hence $f_X$ is strictly decreasing on $\intervalleof{0}{\frac1m}$. 
Finally $f_X(p)\in\intervalleff{1-p}{1-p^2}$ gives continuity at $0$, where $f_X(0)=1>f_X(p)$.

We now prove \ref{it:universal bound Fxp}. Let $p\in\intervalleof{0}{\frac1m}$, so that $p(1-p)m\leq 1-p<1$. 
By convexity of $\phi$, for every $x\in\intervalleff{0}{1}$, we have $1-m (1-x) \leq 	 \phi(x) \leq 1$ so that
	\begin{align*}
1-p+p(1-p)\left(1-m(1-x)\right)-x \leq 1-p+p(1-p)\phi(x)-x \leq 1-p+p(1-p)-x.
	\end{align*}
	Taking $x=f_X(p)$ and using the fixed point equation~\eqref{eq:equation f_X(p)}, we get
	\[
	1-p+p(1-p)\left(1-m(1-f_X(p))\right)-f_X(p)\leq 0\leq 1-p+p(1-p)-f_X(p).
	\]
	The functions $x\mapsto 1-p+p(1-p)\left(1-m(1-x)\right)-x$ and $x\mapsto 1-p+p(1-p)-x$ are affine with negative slopes $p(1-p)m-1$ and $-1$, so that $f_X(p)$ lies between their zeros, which are respectively $1-\frac{p^2}{1-p(1-p)m}$ and $1-p^2$.
	The upper bound is strict: since $f_X(p)\leq 1-p^2<1$, we have $\phi(f_X(p))<1$, so that the fixed point equation~\eqref{eq:equation f_X(p)} gives $f_X(p)<1-p+p(1-p)=1-p^2$.
	Finally, the inequality $p(1-p)m\leq 1-p$ gives $\frac{p^2}{1-p(1-p)m}\leq p$.
	This proves \ref{it:universal bound Fxp}.

	Let us turn to the proof of \ref{it:p,v mapsto Fpv is continuous}. Continuity is clear. A direct computation gives
	\begin{align*}
	\partial_vF_m(p,v)&=(1-p)^2\left(1-mp+2p^2(m-1)v\right),\\
	\partial_pF_m(p,v)&=v(1-p)\left[-(2+m)+3mp+2p(1-2p)(m-1)v\right].
	\end{align*}
	The first quantity is positive for $p<\frac1m$ and $v\geq0$. For $v\in\intervalleof{0}{m}$ and $p\in\intervalleff{0}{\frac1{m^2}}$, the bracket in the second one is at most $-(2+m)+\frac3m+\frac{2(m-1)}{m}=-m+\frac1m<0$. 
\end{proof}

\begin{proof}[Proof of Theorem~\ref{prop:the critical percolation parameter is the solution to an equation}]
We first prove the continuity and the monotonicity.
Since $\phi'(1)=m$ is finite, the function $\phi'$ is continuous on $\intervalleff{0}{1}$, so that $p\mapsto F_m(p,\phi'(f_X(p)))$ is continuous on $\intervalleff{0}{\frac1{m^2}}$ as a composition of continuous functions, by Lemma~\ref{lem:proprietes f_X et F_m}\ref{it:p mapsto fXp is continuous decreasing} and~\ref{it:p,v mapsto Fpv is continuous}.
Moreover, the function $\phi'$ is non-decreasing, and it is positive on $\intervalleof{0}{1}$ because $\Pp{X\geq1}>0$.
Since $f_X$ is decreasing with values in $\intervalleof{0}{1}$ on this interval, the function $p\mapsto\phi'(f_X(p))$ is non-increasing with values in $\intervalleof{0}{m}$.
Hence, for $0\leq p_1<p_2\leq\frac1{m^2}$, Lemma~\ref{lem:proprietes f_X et F_m}\ref{it:p,v mapsto Fpv is continuous} gives
\[
F_m\left(p_2,\phi'(f_X(p_2))\right)\leq F_m\left(p_2,\phi'(f_X(p_1))\right)<F_m\left(p_1,\phi'(f_X(p_1))\right).
\]
We now determine the sign of this function at both ends of the interval.
At $p=0$, we have $f_X(0)=1$, hence $F_m(0,\phi'(f_X(0)))=F_m(0,m)=m-(m-1)=1>0$.
At $p=\frac1{m^2}$, since $\phi'\leq m$ on $\intervalleff{0}{1}$ and $F_m$ is increasing in its second variable, we have
\[
F_m\left(\frac1{m^2},\phi'\left(f_X\left(\frac1{m^2}\right)\right)\right)\leq F_m\left(\frac1{m^2},m\right)=(m-1)\left[\left(1-\frac1{m^2}\right)^2\left(1+\frac1{m^2}\right)-1\right]<0,
\]
where the last inequality comes from the identity $\left(1-\frac1{m^2}\right)^2\left(1+\frac1{m^2}\right)=\left(1-\frac1{m^2}\right)\left(1-\frac1{m^4}\right)$, whose right-hand side is smaller than $1$.
By the intermediate value theorem, the function $p\mapsto F_m(p,\phi'(f_X(p)))$ has a unique zero in $\intervalleff{0}{\frac1{m^2}}$, which lies in $\intervalleoo{0}{\frac1{m^2}}$; moreover, it is positive before this zero and negative after it.

By Lemma~\ref{lem:p is supercritical iff F<0}, for $p\in\intervalleff{0}{\frac1{m^2}}$ there is percolation if and only if $p$ is larger than this zero, so that $p_c(X)$ is equal to this zero.
\end{proof}

\begin{remark}\label{rem:binary tree}
For $X=X^{(2)}$, that is when $T$ is the binary tree, we have $\phi(x)=x^2$, so that the fixed point equation~\eqref{eq:equation f_X(p)} is quadratic and gives
\[
f_{X^{(2)}}(p)=\frac{1-\sqrt{1-4p(1-p)^2}}{2p(1-p)}\qquad\text{for every } p\in\intervalleoo{0}{1}.
\]
Since $\phi'(x)=2x$, Theorem~\ref{prop:the critical percolation parameter is the solution to an equation} shows that $p_c(X^{(2)})$ is the unique solution in $\intervalleoo{0}{\frac14}$ of the explicit equation $F_2\left(p,2f_{X^{(2)}}(p)\right)=0$, which is an algebraic equation in $p$, and we find numerically $p_c(X^{(2)})\approx 0.1586563260$.
\end{remark}

\subsection{Asymptotics of \texorpdfstring{$p_c(X^{(m)})$}{p\_c(X\textasciicircum(m))} for large \texorpdfstring{$m$}{m}}
\begin{proposition}\label{prop:asymptotics pc Xm}
As $m\to\infty$, we have
\[
p_c(X^{(m)})=\frac{1}{m^2}-\frac{1}{m^3}+O\left(\frac1{m^4}\right).
\]
\end{proposition}
\begin{proof}
Throughout the proof, $\phi$ denotes the generating function of $X^{(m)}$ and we assume that $m\geq2$.
Since we expect $p_c(X^{(m)})$ to be of order $\frac1{m^2}$, we evaluate the function $p\mapsto F_m\left(p,\phi'(f_{X^{(m)}}(p))\right)$ of Theorem~\ref{prop:the critical percolation parameter is the solution to an equation} at the point $\frac{\theta}{m^2}$, where $\theta\in\intervalleff{\frac12}{1}$ is a parameter.
All the $O(\cdot)$ terms that appear in the proof below hold as $m\rightarrow \infty$ and uniformly in $\theta\in\intervalleff{\frac12}{1}$.

\step{estimate of $f_{X^{(m)}}$}\label{step:estimate f Xm}
Since $\frac{\theta}{m^2}\left(1-\frac{\theta}{m^2}\right)m\leq\frac1m$, Lemma~\ref{lem:proprietes f_X et F_m}\ref{it:universal bound Fxp} gives
\[
\frac{\theta^2}{m^4}\leq 1-f_{X^{(m)}}\left(\frac{\theta}{m^2}\right)\leq\frac{\theta^2}{m^4}\cdot\frac{1}{1-\frac1m}\leq\frac{\theta^2}{m^4}\left(1+\frac2m\right),
\]
where the last inequality uses $m\geq2$. In other words,
\begin{equation}\label{eq:estimate f Xm}
1-f_{X^{(m)}}\left(\frac{\theta}{m^2}\right)=\frac{\theta^2}{m^4}+O\left(\frac1{m^5}\right).
\end{equation}

\step{estimate of $\phi'\circ f_{X^{(m)}}$}\label{step:estimate phi' f Xm}
Since $X^{(m)}\leq m+1$ almost surely, we have
\begin{align*}
\phi''(1)&=\Ec{X^{(m)}(X^{(m)}-1)}=m^2-m+\{m\}(1-\{m\})=m^2+O(m),\\
\phi'''(1)&=\Ec{X^{(m)}(X^{(m)}-1)(X^{(m)}-2)}\leq (m+1)m(m-1)\leq m^3 .
\end{align*}
Moreover, the function $\phi'$ is convex and $\phi'''$ is non-decreasing on $\intervalleff{0}{1}$, so that Taylor's formula gives, for every $x\in\intervalleff{0}{1}$,
\begin{equation}\label{eq:taylor phi'}
m-\phi''(1)(1-x)\leq\phi'(x)\leq m-\phi''(1)(1-x)+\frac{\phi'''(1)}{2}(1-x)^2 .
\end{equation}
We apply these inequalities at the point $x=f_{X^{(m)}}\left(\frac{\theta}{m^2}\right)$.
By~\eqref{eq:estimate f Xm}, we have $\phi''(1)(1-x)=\frac{\theta^2}{m^2}+O(m^{-3})$ and $\phi'''(1)(1-x)^2=O(m^{-5})$, hence
\begin{equation}\label{eq:estimate phi' f Xm}
\phi'\left(f_{X^{(m)}}\left(\frac{\theta}{m^2}\right)\right)=m-\frac{\theta^2}{m^2}+O\left(\frac1{m^3}\right).
\end{equation}

\step{expansion of $F_m$}\label{step:expansion F Xm}
Recall from~\eqref{eq:definitionF_m} that $F_m(p,v)+(m-1)$ is the product of the three factors $(1-p)^2$, $1-mp+p^2(m-1)v$ and $v$.
When $p$ is replaced by $\frac{\theta}{m^2}$ and $v$ by $\phi'\left(f_{X^{(m)}}\left(\frac{\theta}{m^2}\right)\right)$, the estimate~\eqref{eq:estimate phi' f Xm} shows that these three factors satisfy
\begin{align*}
\left(1-\frac{\theta}{m^2}\right)^2&=1-\frac{2\theta}{m^2}+O\left(\frac1{m^4}\right),\\
1-\frac{\theta}{m}+\frac{\theta^2}{m^4}(m-1)\,\phi'\left(f_{X^{(m)}}\left(\frac{\theta}{m^2}\right)\right)&=1-\frac{\theta}{m}+\frac{\theta^2}{m^2}+O\left(\frac1{m^3}\right),\\
\phi'\left(f_{X^{(m)}}\left(\frac{\theta}{m^2}\right)\right)&=m-\frac{\theta^2}{m^2}+O\left(\frac1{m^3}\right).
\end{align*}
Multiplying these three expansions and subtracting $m-1$, we obtain
\begin{equation}\label{eq:expansion F Xm}
F_m\left(\frac{\theta}{m^2},\phi'\left(f_{X^{(m)}}\left(\frac{\theta}{m^2}\right)\right)\right)=1-\theta+\frac{\theta^2-2\theta}{m}+O\left(\frac1{m^2}\right).
\end{equation}

\step{location of $p_c(X^{(m)})$}\label{step:location pc Xm}
Let $C\geq1$ be a constant such that, for every $m$ large enough and every $\theta\in\intervalleff{\frac12}{1}$, the term $O\left(\frac1{m^2}\right)$ in~\eqref{eq:expansion F Xm} is at most $\frac{C}{m^2}$ in absolute value, and assume from now on that $m$ is large enough for this bound to hold, and that $m\geq\max(2C,4)$.
Since $\theta^2-2\theta=(1-\theta)^2-1$, the right-hand side of~\eqref{eq:expansion F Xm} is equal to
\begin{equation}\label{eq:expansion F Xm bis}
(1-\theta)-\frac1m+\frac{(1-\theta)^2}{m}+O\left(\frac1{m^2}\right).
\end{equation}
We use this expression for two values of $\theta$, which both belong to $\intervalleff{\frac12}{1}$ thanks to our assumption on $m$.
\begin{itemize}
\item For $\theta=1-\frac1m+\frac{2C}{m^2}$, we have $0\leq 1-\theta\leq\frac1m$, so that~\eqref{eq:expansion F Xm bis} is at most $-\frac{2C}{m^2}+\frac1{m^3}+\frac{C}{m^2}$, which is negative since $C\geq1$.
\item For $\theta=1-\frac1m-\frac{2C}{m^2}$, the quantity~\eqref{eq:expansion F Xm bis} is at least $\frac{2C}{m^2}-\frac{C}{m^2}$, which is positive.
\end{itemize}
By Theorem~\ref{prop:the critical percolation parameter is the solution to an equation}, the function $p\mapsto F_m\left(p,\phi'(f_{X^{(m)}}(p))\right)$ is strictly decreasing on $\intervalleff{0}{\frac1{m^2}}$ and vanishes at $p_c(X^{(m)})$. 
Hence
\[
\frac{1}{m^2}\left(1-\frac1m-\frac{2C}{m^2}\right)<p_c(X^{(m)})<\frac{1}{m^2}\left(1-\frac1m+\frac{2C}{m^2}\right),
\]
which proves the proposition.
\end{proof}

\section{Monotonicity in the equation}\label{sec:monotonicity}
Recall that $X$ is  a random variable with non-negative integer values that has expectation $m>1$ and that $\phi(x)=\Ec{x^X}$ is its generating function.

\subsection{Some criteria for the comparison of critical parameters}
Now, consider $Y$ another random variable with non-negative integer values with expectation $m>1$.
We denote by $\psi(x)=\Ec{x^Y}$ its generating function.
Then, if we denote by $p_c(X)$ and $p_c(Y)$ the respective associated critical parameters, the results of the previous section entail the following result.

\begin{proposition}\label{prop:criterion for comparison of critical points}
Denoting by $\mathrm{sgn}(t)\in\{-1,0,1\}$ the sign of a real number $t$, we have
	\begin{align*}
	\mathrm{sgn}\left(p_c(X)-p_c(Y)\right)=\mathrm{sgn}\left(\phi'\left(f_X(p_c(X))\right)-\psi'\left(f_Y(p_c(X))\right)\right).
	\end{align*}
\end{proposition}
\begin{proof}
Recall that $X$ and $Y$ both have mean $m$. 
By Theorem~\ref{prop:the critical percolation parameter is the solution to an equation} applied to $Y$, the function $p\mapsto F_m\left(p,\psi'(f_Y(p))\right)$ is strictly decreasing on $\intervalleff{0}{\frac1{m^2}}$ and vanishes at $p_c(Y)$. 
Since $p_c(X)$ also belongs to the interval $\intervalleff{0}{\frac1{m^2}}$, we get
\[
\mathrm{sgn}\left(p_c(Y)-p_c(X)\right)=\mathrm{sgn}\left(F_m\left(p_c(X),\psi'(f_Y(p_c(X)))\right)\right).
\]
On the other hand, the function $F_m(p_c(X),\cdot)$ is strictly increasing on $\intervallefo{0}{\infty}$ by Lemma~\ref{lem:proprietes f_X et F_m}\ref{it:p,v mapsto Fpv is continuous}, and it vanishes at $\phi'(f_X(p_c(X)))$ by Theorem~\ref{prop:the critical percolation parameter is the solution to an equation} applied to $X$. 
Hence
\[
\mathrm{sgn}\left(F_m\left(p_c(X),\psi'(f_Y(p_c(X)))\right)\right)=\mathrm{sgn}\left(\psi'\left(f_Y(p_c(X))\right)-\phi'\left(f_X(p_c(X))\right)\right).
\]
Combining these two identities gives the result.
\end{proof}
The rest of the paper relies on the following corollary of the above proposition, using the fact that $p_c(X)\in\intervalleoo{0}{\frac1{m^2}}$ by Theorem~\ref{prop:the critical percolation parameter is the solution to an equation}.
\begin{corollary}\label{cor:sufficient condition for comparison}
	If for all $p\in \intervalleoo{0}{\frac{1}{m^2}}$ we have
	\begin{align*}
		\psi'(f_Y(p)) < \phi'(f_X(p)),
	\end{align*}
then $p_c(X) >  p_c(Y)$.
\end{corollary}

A first natural heuristic in order to verify the above condition is to look at what happens when $p$ is small.
Assume that $X$ has a finite second moment.
In that case we can Taylor expand $f_X(\cdot)$ around $p=0$ and $\phi'(\cdot)$ around $1$, using the estimate $1-f_X(p)=p^2+O(p^3)$ given by Lemma~\ref{lem:proprietes f_X et F_m}\ref{it:universal bound Fxp}, and get that
\begin{align}\label{eq:taylor expansion phi prime C p phi}
	\phi'(f_X(p)) \underset{p \rightarrow 0}{=} m - \phi''(1)p^2 + o(p^2).
\end{align}
Since $\phi''(1) = \Ec{X(X-1)}= \Var{X} + m(m-1)$, this entails that the sufficient condition given in Corollary~\ref{cor:sufficient condition for comparison} can only hold if $\Var{X} \leq \Var{Y}$.
This suggests that somehow, the random variables $X$ and $Y$ for which we will achieve some comparison using Corollary~\ref{cor:sufficient condition for comparison} satisfy that $Y$ is more variable than $X$ in some sense.

\subsection{The convex order case}
The previous subsection leads us to consider the case where $X$ and $Y$ are comparable in the convex order, defined in Section~\ref{sec:intro model}.
In this subsection, we always assume that $X$ and $Y$ have mean $m>1$, and that $X\preceq_{\mathrm{cx}} Y$.

Recall the definition \eqref{eq:def loi X^{(m)}} of $X^{(m)}$. 
The following lemma ensures that $X^{(m)}$ is minimal for the convex order, among random variables with mean $m$. 
We provide a proof for completeness. 
\begin{lemma}\label{lem:Xm minimal}
For every $\Z_+$-valued random variable $Y$ with mean $m$, we have $X^{(m)}\preceq_{\mathrm{cx}}Y$.
\end{lemma}
\begin{proof}
Let $f:\Z_+\to\R$ be convex and let $L$ be the affine function which coincides with $f$ at $\lfloor m\rfloor$ and at $\lfloor m\rfloor+1$. 
By convexity, $f(j)\geq L(j)$ for every $j\in\Z_+$. 
Hence $\Ec{f(Y)}\geq\Ec{L(Y)}=L(m)=(1-\{m\})f(\lfloor m\rfloor)+\{m\}f(\lfloor m\rfloor+1)=\Ec{f(X^{(m)})}$.
\end{proof}

First, we remark that for any fixed $x\in \intervalleff{0}{1}$, the function $k\mapsto x^k$ is convex so that
\begin{align}\label{eq:phi leq psi}
\phi(x) = \Ec{x^X} \leq \Ec{x^Y} = \psi(x) .
\end{align}
Plugging inequality~\eqref{eq:phi leq psi} into the fixed point equation \eqref{eq:equation f_X(p)} we get that for any $p$, we have
\begin{align}\label{eq:inequality C p phi}
f_X(p) \leq f_Y(p).
\end{align}
The next lemma, which is used in the proof of Lemma~\ref{lem:if N small then comparison} below, states that inequality~\eqref{eq:phi leq psi} is strict on $\intervalleoo{0}{1}$ when $X\prec_{\mathrm{cx}}Y$.
\begin{lemma}\label{lem:strict phi psi}
Assume that $X\prec_{\mathrm{cx}}Y$. 
Then $\phi(x)<\psi(x)$ for every $x\in\intervalleoo{0}{1}$.
\end{lemma}
\begin{proof}
Let $x\in\intervalleoo{0}{1}$. 
For any integer $k\geq 0$, writing $x^k -1$ as the telescopic sum $\sum_{i=0}^{k-1} x^{i} (x-1)$ we have
\begin{align*}
x^k-1-k(x-1)=(x-1)\sum_{i=0}^{k-1}\left(x^i-1\right)=(x-1)^2\sum_{i=0}^{k-1}\sum_{j=0}^{i-1}x^j=(1-x)^2\sum_{j=0}^{\infty}x^j\,(k-j-1)_+,
\end{align*}
where the last equality comes from the fact that, for every $j\geq0$, the number of integers $i$ such that $j<i\leq k-1$ is $(k-j-1)_+$.
We apply this identity with $k=X$ and with $k=Y$, and take expectations; all the terms of the sum are non-negative, so that we can exchange the sum and the expectation.
Since $X$ and $Y$ have the same mean $m$, we get
\[
\psi(x)- \phi(x)= (1-x)^2\sum_{j=0}^\infty x^j \left(\Ec{(Y-j-1)_+} - \Ec{(X-j-1)_+}\right).
\]
Every term of this sum is non-negative by the convex ordering $X\preceq_{\mathrm{cx}}Y$, since $k\mapsto (k-j-1)_+$ is convex. 
If for some $x\in \intervalleoo{0}{1}$ we had  $\psi(x)=\phi(x)$, then it would mean that all these terms vanish, so that $\Ec{(Y-j)_+}=\Ec{(X-j)_+}$ for every $j\geq1$
and also for $j=0$ since $X$ and $Y$ have the same mean. 
Since $\Pp{X>j}=\Ec{(X-j)_+}-\Ec{(X-j-1)_+}$ for every $j\geq 0$, and similarly for $Y$, the random variables $X$ and $Y$ would then have the same law, which is excluded.
\end{proof}

Now for $p\in \intervalleof{0}{\frac{1}{m^2}}$, we write
\begin{align}\label{eq:writing as good term plus bad term}
	\psi'(f_Y(p))- \phi'(f_X(p))  = ( \underset{\text{good term}}{\underbrace{\psi'(f_Y(p)) - \phi'(f_Y(p))}}
	) 	+( \underset{\text{bad term}}{\underbrace{\phi'(f_Y(p)) - \phi'(f_X(p))}})
\end{align}
and try to bound the two terms on the right-hand side separately.
Since we aim at proving that the quantity on the left-hand side is negative, the second term is ``bad'' in the sense that it is non-negative (because of \eqref{eq:inequality C p phi} and the monotonicity of $\phi'$), and we count on the first one to be negative enough to make up for it.

\paragraph{Outline of proof of Theorem~\ref{thm:convex order + bound support + large enough mean = comparison}.}
We want to compare the contributions from the good and the bad term.
The trick here will be to upper-bound the two quantities in such a way that we can factor out the non-negative term $\psi(f_Y(p))-\phi(f_Y(p))$.
For the bad term, this follows from the fixed point equation~\eqref{eq:equation f_X(p)} and holds for every $p<\frac1m$, see Lemma~\ref{lem:controlling the bad term}.
For the good term, this is where we use the bound on the support of $X$, see Lemma~\ref{lem:inequality difference derivatives}.
We are then left with determining the sign of an explicit factor, which is done in Lemma~\ref{lem:if N small then comparison} and Lemma~\ref{lem:bound Nm}.

\subsubsection{Controlling the bad term}

\begin{lemma}\label{lem:controlling the bad term}
Assume that $\phi''(1)<\infty$. For $p < 1/m$, we have
\begin{align}\label{eq:upperbound bad term}
	\phi'(f_Y(p)) - \phi'(f_X(p))	\leq \phi''(1) \cdot \frac{p(1-p)}{1-p(1-p)m} \cdot  \left( \psi(f_Y(p)) - \phi(f_Y(p)) \right).
\end{align}
\end{lemma}
\begin{proof}
Using the mean value theorem, inequality~\eqref{eq:inequality C p phi} and the fact that $\phi''$ is non-decreasing on $\intervalleff{0}{1}$, we have
\begin{align}\label{eq:TAF bad term}
	\phi'(f_Y(p)) - \phi'(f_X(p)) \leq \phi''(1) \cdot (f_Y(p) - f_X(p)).
\end{align}
Now we work on bounding the right-hand side.
Subtracting the two fixed-point equations \eqref{eq:equation f_X(p)} satisfied respectively by $f_X(p)$ and $f_Y(p)$, corresponding to respectively $\phi$ and $\psi$, yields
\begin{align*}
	f_Y(p) - f_X(p) &= p(1-p) \left[ \psi(f_Y(p)) - \phi(f_X(p)) \right]\\
	&= p(1-p) \left[ \psi(f_Y(p)) - \phi(f_Y(p)) + \phi(f_Y(p)) - \phi(f_X(p)) \right] \\
	&\leq p(1-p) \left[ \psi(f_Y(p)) - \phi(f_Y(p)) \right] + p(1-p)m(f_Y(p) - f_X(p)),
\end{align*}
where the inequality comes from the mean value theorem, inequality~\eqref{eq:inequality C p phi} and the fact that $\phi'\leq m$ on $\intervalleff{0}{1}$.
This yields
\begin{align*}
 f_Y(p)- f_X(p) \leq \frac{p (1-p)}{1-p(1-p)m }\cdot \left(\psi(f_Y(p)) - \phi(f_Y(p))\right).
\end{align*}
Plugging this back into \eqref{eq:TAF bad term}, we get~\eqref{eq:upperbound bad term}.
\end{proof}

\subsubsection{Controlling the good term}

Now, we focus on the ``good term'' $\psi'(f_Y(p)) - \phi'(f_Y(p))$, which is expressed as the difference of the two derivatives evaluated at the same value $f_Y(p)$.
This motivates the next lemma, in which we crucially assume a bound on the support of $X$.

\begin{lemma}\label{lem:inequality difference derivatives}
	Suppose that $X\preceq_{\mathrm{cx}} Y$ and that $\supp X \subset \{0,\dots ,N\}$ for some integer $N\geq1$.
	Then for all $x\in\intervallefo{\frac12}{1}$, we have
	\begin{align*}
		\psi'(x) - \phi'(x) \leq \frac{1}{x}\left(\frac{1}{1-x} - N\right) \cdot(\phi(x) - \psi(x)).
	\end{align*}
\end{lemma}
\begin{proof}
	We introduce the function $g:\intervallefo{0}{\infty}\to\R$ defined by
	\begin{align*}
		g(s)= (s+a) x^{s-1},
	\end{align*}
	where $a=a(N,x)$ is a real parameter that we will choose later on.
	Note that $\Ec{g(Y)} = \psi'(x) + a \frac{\psi(x)}{x}$.
	We compute the first and second derivative of $g$ and get
	\begin{align*}
		g'(s)= x^{s-1}(1+(s+a) \log x)  \quad \text{and} \quad 	g''(s)= x^{s-1} (\log x) (2+(s+a) \log x).
	\end{align*}
	Now we choose $a$ in such a way that $g(N-1)=g(N)$\footnote{An earlier version of the argument chose $a$ so that the function $g$ is maximized at $N$. This trick was suggested to us by Gemini 3.1 Pro.}, that is $(N-1+a) x^{N-2} = (N+a)x^{N-1}$.
	This gives $N-1+a = Nx + ax$ and so
	\begin{align*}
		a= \frac{N(1-x)-1}{x-1}= \frac{1}{1-x} - N.
	\end{align*}
	Since the sign of $g'(s)$ is that of $1+(s+a)\log x$, which is decreasing in $s$, the function $g$ is increasing then decreasing, and $g(N-1)=g(N)$ implies that its maximum is reached in $\intervalleff{N-1}{N}$; in particular $g$ is decreasing on $\intervallefo{N}{\infty}$. 
	Let us check that $g$ is concave on $\intervalleff{0}{N}$.
	We have $g''(s)\leq 0$ if and only if $s+a\leq\frac{2}{-\log x}$. For $s\leq N$ we have $s+a\leq N+a=\frac{1}{1-x}$, and $\frac1{1-x}\leq\frac{2}{-\log x}$ because the function $x\mapsto 2(1-x)+\log x$ vanishes at $1$ and has derivative $\frac1x-2\leq0$ on $\intervalleff{\frac12}{1}$. 
	Hence $g$ is concave on $\intervalleff{0}{N}$. 
	This is where we use $x\geq\frac12$; the argument works as long as $-\log x\leq 2(1-x)$, that is $x\gtrsim 0.2032$.
	We introduce
	the function $g^*$ defined on the non-negative integers as
	\begin{equation*}
		g^*(k) =
		\left \lbrace
		\begin{aligned}
			&g(k) \quad \text{if} \quad k \leq N\\
			&g(N) \quad \text{otherwise}.
		\end{aligned}
		\right.
	\end{equation*}
	Note that $g^*$ is concave and $g\leq g^*$: the increments $g^*(k+1)-g^*(k)$ are non-increasing for $k\leq N-1$ by concavity of $g$ on $\intervalleff{0}{N}$, they vanish for $k\geq N-1$ since $g(N-1)=g(N)$, and $g(k)\leq g(N)$ for $k\geq N$ since $g$ is decreasing on $\intervallefo{N}{\infty}$. 
	Note also that $g^*$ is bounded.
	Using in turn the inequality $g\leq g^*$, the convex ordering $X\preceq_{\mathrm{cx}} Y$ applied to the convex function $-g^*$, and then the fact that $g^*(X)=g(X)$ a.s., we get
	\begin{align*}
		\psi'(x) + a \frac{\psi(x)}{x}= \Ec{g(Y)}\leq \Ec{g^*(Y)} \leq \Ec{g^*(X)}= \Ec{g(X)} =  \phi'(x) + a \frac{\phi(x)}{x}
	\end{align*}
	so that
	\begin{align}\label{eq:quantitative inequality on derivatives}
		\psi'(x) - \phi'(x) \leq  \frac{a}{x} \cdot(\phi(x) - \psi(x)),
	\end{align}
	which is the result, since $\frac ax=\frac1x\left(\frac1{1-x}-N\right)$.
\end{proof}

\subsubsection{Summing good term and bad term}
For any $p\in\intervallefo{0}{\frac1m}$ and $C\in \intervalleoo{0}{1}$, let
\begin{align}\label{eq:definition Nm}
	N_m(p,C) := \frac{1- p(1-p) m \left( 1 - C + C^2 \right) }{\left( 1-C \right) \left( 1-p(1-p)m(1-C) \right)}.
\end{align}
\begin{lemma}
	\label{lem:if N small then comparison}
	Suppose that $\supp X \subset \{0,1,\dots ,N\}$, that $m\geq\sqrt2$, and that $X\prec_{\mathrm{cx}}Y$.
For any $p\in\intervalleof{0}{\frac1{m^2}}$, if $N < N_m(p,f_Y(p))$ then $\psi'(f_Y(p))< \phi'(f_X(p))$.
\end{lemma}
\begin{proof}
First, Lemma~\ref{lem:proprietes f_X et F_m}\ref{it:universal bound Fxp} gives $f_Y(p)\geq 1-p\geq 1-\frac1{m^2}\geq\frac12$, so that we can apply Lemma~\ref{lem:inequality difference derivatives} with $x=f_Y(p)$.
We also apply Lemma~\ref{lem:controlling the bad term}, whose assumption $\phi''(1)<\infty$ holds since $X$ has bounded support.
Plugging these two bounds into the decomposition~\eqref{eq:writing as good term plus bad term} and factoring out the term $\psi(f_Y(p))-\phi(f_Y(p))$, we get
\begin{align} \label{eq:good term plus bad term}
	\psi'(f_Y(p))- \phi'(f_X(p)) \leq \left(\phi''(1)\frac{p(1-p)}{1-p(1-p)m}  -\frac{1}{f_Y(p)}\left(\frac{1}{1-f_Y(p)}-N\right) \right) \cdot \left( \psi(f_Y(p))-\phi(f_Y(p))\right).
\end{align}
Note that since $f_Y(p) \in \intervalleoo{0}{1}$ the term $\left( \psi(f_Y(p))-\phi(f_Y(p))\right)$ is positive by Lemma~\ref{lem:strict phi psi}.
The sign of the right-hand side of \eqref{eq:good term plus bad term} is hence determined by the sign of its first factor, which satisfies
\begin{align} \label{eq:firstbound}
	\phi''(1)\frac{p(1-p)}{1-p(1-p)m}  -\frac{1}{f_Y(p)}\left(\frac{1}{1-f_Y(p)}-N\right) 
	& \leq \frac{p(1-p) m (N-1)}{1-p(1-p)m} - \frac{1}{f_Y(p)} \left( \frac{1}{1-f_Y(p)} -N \right),
\end{align}
where we used
\[
\phi''(1) = \sum_{k=0}^{N} k (k-1) \mu(k) \leq (N-1) \sum_{k=0}^{N} k \mu(k) = m(N-1).
\]

We can write the right-hand side of our upper-bound \eqref{eq:firstbound} in the product form
\begin{align}\label{eq:firstbound product form}
 \frac{p(1-p) m (N-1)}{1-p(1-p)m} - \frac{1}{f_Y(p)} \left( \frac{1}{1-f_Y(p)} -N \right) &=
 \frac{1-p(1-p)m(1-f_Y(p))}{(1-p(1-p)m) \, f_Y(p)} \left( N - N_m(p,f_Y(p)) \right)
\end{align}
where $N_m(p,C)$ is defined in \eqref{eq:definition Nm}.
Since
\begin{align*}
	 \frac{1-p(1-p)m(1-f_Y(p))}{(1-p(1-p)m) \, f_Y(p)}>0,
\end{align*}
the sign of the whole expression \eqref{eq:firstbound product form} is that of the difference $\left( N - N_m(p,f_Y(p)) \right)$.
In particular, if $N < N_m(p,f_Y(p))$ then the expression is negative and hence $\psi'(f_Y(p))- \phi'(f_X(p))< 0$.
\end{proof}
In order to know how large $N$ can be, we now just need to provide a lower-bound on  $N_m(p,f_Y(p))$.
This is done in the next lemma.
\begin{lemma}\label{lem:bound Nm}
For any $p\in \intervalleof{0}{\frac{1}{m^2}}$ we have
\begin{align*}
		N_m(p,f_Y(p)) > 1+\frac{(m-1)^2(m^3+1)}{m}.
\end{align*}
\end{lemma}
\begin{proof}
A simple Euclidean division yields
\begin{align*}
	N_m(p,C) = 1 + \frac{(1-p(1-p) m)  C }{\left( 1-C \right) \left( 1-p(1-p)m(1-C) \right)}.
\end{align*}
Using $0<p\leq \frac{1}{m^2}$ and $C\in \intervalleoo{0}{1}$, we get 
\begin{align*}
	N_m(p,C)> 1 + (1-p(1-p) m) \cdot C \cdot \frac{1 }{ 1-C}.
\end{align*}

By Lemma~\ref{lem:proprietes f_X et F_m}\ref{it:universal bound Fxp}, we have $p^2\leq 1-f_Y(p)\leq\frac{p^2}{1-p(1-p)m}$, so that
\[
f_Y(p)\geq\frac{1-p(1-p)m-p^2}{1-p(1-p)m} =\frac{1-pm + p^{2}(m-1)}{1-p(1-p)m} >0\qquad \text{and} \qquad \frac{1}{1-f_Y(p)}\geq\frac{1-p(1-p)m}{p^2}>0.
\]
Multiplying these two lower bounds of positive quantities, we get
\begin{align*}
N_m(p,f_Y(p))
&> 1+\left(1-p(1-p)m\right)\cdot\frac{1-p(1-p)m-p^2}{1-p(1-p)m}\cdot\frac{1-p(1-p)m}{p^2}\\&
=1+\left(1-p(1-p)m\right)\,\frac{1-p(1-p)m-p^2}{p^2}.
\end{align*}
Now we bound the two factors separately. 
For the first factor, we just write $1-p(1-p)m\geq 1-\frac1m=\frac{m-1}{m}$. 
For the second factor, we write
\begin{align*}
	\frac{1-p(1-p)m-p^2}{p^2}=\frac{1}{p^2}-\frac mp+m-1.
\end{align*}
Its derivative in $p$ is  $\frac{mp-2}{p^3}$, which is negative for $p\in \intervalleoo{0}{\frac{1}{m^2}}$, so we bound this factor from below by its value at $p=\frac{1}{m^2}$, which is $m^4-m^3+m-1=(m-1)(m^3+1)$. This gives the result.
\end{proof}

Now we just need to put everything together to obtain a proof of our main result.
\begin{proof}[Proof of Theorem~\ref{thm:convex order + bound support + large enough mean = comparison}]
Let $m>1$. The function $m\mapsto \frac{(m-1)^2(m^3+1)}{m}=(m-1)^2\left(m^2+\frac1m\right)$ is increasing on $\intervallefo{1}{\infty}$, hence $N(m)\geq2$ if and only if $\frac{(m-1)^2(m^3+1)}{m}\geq 1$, that is $m\geq m_1$. If $m<m_1$, there is no random variable with mean $m>1$ supported by $\{0,\dots,N(m)\}\subset\{0,1\}$, and there is nothing to prove.

Assume now that $m\geq m_1$, that $X\prec_{\mathrm{cx}}Y$ and that $\supp X \subset \{0,1,\dots ,N(m)\}$. 
Note that $m_1>\sqrt2$, so that the assumption $m\geq\sqrt2$ of Lemma~\ref{lem:if N small then comparison} is satisfied. Thanks to Lemma~\ref{lem:bound Nm} and Lemma~\ref{lem:if N small then comparison}, for any $p\in\intervalleof{0}{\frac{1}{m^2}}$ we have $N(m)\leq 1+\frac{(m-1)^2(m^3+1)}{m}< N_m(p,f_Y(p))$ and thus $\psi'(f_Y(p))< \phi'(f_X(p))$.
Thanks to Corollary~\ref{cor:sufficient condition for comparison}, we can conclude that $p_c(X)> p_c(Y)$.
\end{proof}
\begin{proof}[Proof of Corollary~\ref{cor:maximal percolation threseshold}]
Let $m\geq m_1$ and let $Y$ have mean $m$ and a law different from that of $X^{(m)}$.
By Lemma~\ref{lem:Xm minimal}, we have $X^{(m)}\prec_{\mathrm{cx}}Y$, and $\supp X^{(m)}\subset\{0,\dots,\lceil m\rceil\}$.
Let us check that $N(m)\geq\lceil m\rceil$, so that $\supp X^{(m)}\subset\{0,\dots,N(m)\}$.
If $m\in\intervalleff{m_1}{2}$, then we have $N(m)\geq2=\lceil m\rceil$, as seen in the proof of Theorem~\ref{thm:convex order + bound support + large enough mean = comparison}.
If $m\geq 2$, then $(m-1)^2\geq1$ gives $1+\frac{(m-1)^2(m^3+1)}{m}\geq 1+\frac{m^3+1}{m}>m+1$, hence $N(m)\geq\lfloor m+1\rfloor\geq\lceil m\rceil$.
Hence Theorem~\ref{thm:convex order + bound support + large enough mean = comparison} gives $p_c(X^{(m)})> p_c(Y)$.
The asymptotics of $p_c(X^{(m)})$ are given by Proposition~\ref{prop:asymptotics pc Xm}.
\end{proof}
\section{Counter-examples and small critical parameters}\label{sec:counterexamples}
Fix $m>1$.
\subsection{Small critical parameters}
The lower bounds on $p_c$ for trees without leaves, in Theorem~\ref{thm:contre exemple pc tend vers 0}\ref{it:pc without leaves m<2} and in Theorem~\ref{thm:contre-exemple m proche de 1}, rely on the following lemma.
\begin{lemma}\label{lem:lower bound without leaves}
Let $X$ be a random variable with mean $m$ such that $X\geq 1$ almost surely, and denote by $\phi$ its generating function.
For every $p\in\intervalleoo{0}{\frac1m}$, we have
\begin{equation}\label{eq:lower bound F without leaves}
F_m\left(p,\phi'(f_X(p))\right)\geq F_m\left(p,\left(1-\frac{p^2}{1-p(1-p)m}\right)^{m-1}\right).
\end{equation}
\end{lemma}
\begin{proof}
Let $p\in\intervalleoo{0}{\frac1m}$.
Since $s\mapsto f_X(p)^{s-1}$ is convex and $X\geq 1$, Jensen's inequality and Lemma~\ref{lem:proprietes f_X et F_m}\ref{it:universal bound Fxp} give
\[
\phi'(f_X(p))=\Ec{X f_X(p)^{X-1}}\geq\Ec{f_X(p)^{X-1}}\geq f_X(p)^{m-1}\geq\left(1-\frac{p^2}{1-p(1-p)m}\right)^{m-1},
\]
where the last quantity is positive since $p(1-p)m+p^2<1-p+p^2\leq1$.
The result follows since $F_m(p,\cdot)$ is increasing by Lemma~\ref{lem:proprietes f_X et F_m}\ref{it:p,v mapsto Fpv is continuous}.
\end{proof}

Fix an integer $k> m$ and consider $X_k^{(m)}$ supported on $\{1,k\}$ with mean $m$,
 and similarly $Y_k^{(m)}$ supported on $\{0,k\}$ with mean $m$.
We call $\phi_k$ and $\psi_k$ their respective generating functions, so that for all $x\in \intervalleff{0}{1}$
\[
		\phi_k(x) = \frac{k-m}{k-1} x + \frac{m-1}{k-1} x^k \qquad \text{and} \qquad \psi_k(x) =\frac{k-m}{k} +  \frac{m}{k} x^k.
\]

\begin{proof}[Proof of Theorem~\ref{thm:contre exemple pc tend vers 0}]	
\emph{The limit $k\to\infty$.} 
Fix $p\in\intervalleoo{0}{\frac1m}$. 
Since $f_{X_k^{(m)}}(p)$ and $f_{Y_k^{(m)}}(p)$ are at most $1-p^2$ by Lemma~\ref{lem:proprietes f_X et F_m}\ref{it:universal bound Fxp}, and since $k(1-p^2)^{k-1}\to0$ as $k\to\infty$, we get
\begin{align*}
	\phi_k'\left(f_{X_k^{(m)}}(p)\right) = 1-\frac{m-1}{k-1}+\frac{k(m-1)}{k-1}\left(f_{X_k^{(m)}}(p) \right)^{k-1} \xrightarrow[k\to\infty]{}1
\end{align*}
and similarly 
\begin{align*}
	\psi_k'\left(f_{Y_k^{(m)}}(p)\right) = m\left(f_{Y_k^{(m)}}(p) \right)^{k-1} \xrightarrow[k\to\infty]{}0.
\end{align*}
Since the function $F_m(p,\cdot)$ is continuous by Lemma~\ref{lem:proprietes f_X et F_m}\ref{it:p,v mapsto Fpv is continuous}, we deduce that
\[
F_m\left(p,\phi_k'(f_{X_k^{(m)}}(p))\right)\xrightarrow[k\to\infty]{}F_m(p,1)
\qquad\text{and}\qquad
F_m\left(p,\psi_k'(f_{Y_k^{(m)}}(p))\right)\xrightarrow[k\to\infty]{}F_m(p,0).
\]
By Lemma~\ref{lem:p is supercritical iff F<0}, there is percolation at $p$ as soon as the corresponding value of $F_m$ is negative. 

Now, for $p\in\intervalleoo{0}{\frac1m}$, we have $F_m(p,0) = -(m-1) <0$. 
Hence, for $k$ large enough, we have $F_m\left(p,\psi_k'(f_{Y_k^{(m)}}(p))\right)<0$, and therefore $p_c(Y_k^{(m)}) \leq p$.
This proves that $p_c(Y_k^{(m)})\to0$ as $k\to\infty$, so that \ref{it:pc to 0 with leaves} holds.

Similarly for $p\in\intervalleoo{0}{\frac1m}$, we have 
\[
F_m(p,1)=(1-p)^3\left(1-(m-1)p\right)-(m-1)=2-m-\left(1-(1-p)^3\left(1-(m-1)p\right)\right).
\]
For $p\in\intervalleoo{0}{\frac1m}$, both $1-p$ and $1-(m-1)p$ belong to $\intervalleoo{0}{1}$, so that the last term is negative and $F_m(p,1)<2-m$. 

If $m\geq2$, we get $F_m(p,1)<0$ for every $p\in\intervalleoo{0}{\frac1m}$, so that $p_c(X_k^{(m)})\to0$ as $k\to\infty$. 
This proves~\ref{it:pc to 0 without leaves}.

If $1<m<2$, let 
\begin{align*}
	p^+_m := \inf\enstq{p\in\intervalleoo{0}{\frac{1}{m}}}{F_m(p,1)<0}.
\end{align*}
Then the above argument shows that $p_c(X_k^{(m)})\leq p$ for every $p\in\intervalleoo{0}{\frac1m}$ such that $F_m(p,1)<0$ and every $k$ large enough, so that
\begin{align*}
\limsup_{k\to\infty}p_c(X_k^{(m)})\leq p_m^+.
\end{align*}
Now, uniformly in $m\in \intervalleff{1}{2}$ and $p\in \intervalleff{0}{\frac12}$ we have, as $p\rightarrow 0$,
\begin{align*}
	F_m(p,1)=2-m-\left(1-(1-p)^3\left(1-(m-1)p\right)\right) = 2-m -p(m+2) + O(p^2).
\end{align*}
Fix $\epsilon >0$ and plug in $p= \frac{2-m}{m+2} (1+\epsilon)$ to get 
\begin{align*}
	F_m\left(\frac{2-m}{m+2} (1+\epsilon),1\right) =  -(2-m) \epsilon + O((2-m)^2),
\end{align*}
which is negative when $m$ is close enough to $2$.
This entails that, as $m\uparrow 2$ we have
\begin{align*}
p_m^+ \leq \frac{2-m}{m+2} (1+o(1)).
\end{align*}

\emph{Lower bound in~\ref{it:pc without leaves m<2}.}
Let $1<m<2$. We let
\begin{align*}
	p_m^- := \inf\enstq{p\in\intervalleoo{0}{\frac1m}}{F_m\left(p,\left(1-\frac{p^2}{1-p(1-p)m}\right)^{m-1}\right)\leq0}.
\end{align*}
The set in this definition is non-empty: the second argument of $F_m$ in it is at most $1\leq m$, and $F_m$ is increasing in its second variable, so that it contains $\frac1{m^2}$ since $F_m\left(\frac1{m^2},m\right)<0$, as seen in the proof of Theorem~\ref{prop:the critical percolation parameter is the solution to an equation}. 
Note that since $F_m(0,1)= 2-m >0$, by continuity of $p \mapsto F_m\left(p,\left(1-\frac{p^2}{1-p(1-p)m}\right)^{m-1}\right)$ the value $p_m^-$ is positive. 
Now let $X\geq1$ with mean $m$, and let $p\in\intervalleoo{0}{p_m^-}$.
The right-hand side of~\eqref{eq:lower bound F without leaves} is positive, so that there is no percolation at $p$ by Lemma~\ref{lem:lower bound without leaves} and Lemma~\ref{lem:p is supercritical iff F<0}.
Hence $p_c(X) \geq p_m^- >0$. This already shows that
\begin{align*}
\inf\enstq{p_c(X)}{X\geq 1 \text{ a.s.},\ \Ec{X}=m}>0.
\end{align*}

Let us now study $p_m^-$ as $m\uparrow2$. 
For $m\in\intervalleff{1}{2}$ and $p\in\intervalleff{0}{\frac12}$, we have $1-p(1-p)m\geq\frac12$, so that $0\leq\frac{p^2}{1-p(1-p)m}\leq 2p^2\leq\frac12$ and, since $0\leq m-1\leq 1$,
\[
1\geq\left(1-\frac{p^2}{1-p(1-p)m}\right)^{m-1}\geq 1-\frac{p^2}{1-p(1-p)m}\geq 1-2p^2.
\]
Moreover, for $v\in\intervalleff{0}{1}$ we have $0\leq\partial_vF_m(p,v)=(1-p)^2\left(1-mp+2p^2(m-1)v\right)\leq\frac32$, so that the right-hand side of~\eqref{eq:lower bound F without leaves} differs from $F_m(p,1)$ by at most $3p^2$. Using the expansion of $F_m(p,1)$ obtained above, we get that, uniformly in $m\in\intervalleff{1}{2}$ and $p\in\intervalleff{0}{\frac12}$, as $p\to0$,
\begin{align*}
F_m\left(p,\left(1-\frac{p^2}{1-p(1-p)m}\right)^{m-1}\right)=2-m-p(m+2)+O(p^2).
\end{align*}
Fix $\epsilon\in\intervalleoo{0}{1}$. 
For every $p\in\intervalleof{0}{(1-\epsilon)\frac{2-m}{m+2}}$, we have $p(m+2)\leq(1-\epsilon)(2-m)$ and $p\leq 2-m$, so that the left-hand side of the last display is at least $\epsilon(2-m)+O\left((2-m)^2\right)$, uniformly in such $p$. 
It is therefore positive for every such $p$ when $m$ is close enough to $2$, so that $p_m^-\geq(1-\epsilon)\frac{2-m}{m+2}$. 
This entails that, as $m\uparrow2$, we have
\begin{align*}
p_m^-\geq\frac{2-m}{m+2}\left(1-o(1)\right).
\end{align*}

\emph{Conclusion of the proof of~\ref{it:pc without leaves m<2}.} 
Let $1<m<2$. 
Since $X_k^{(m)}\geq1$ a.s.\ and $\Ec{X_k^{(m)}}=m$, the infimum in~\ref{it:pc without leaves m<2} is at most $p_c(X_k^{(m)})$ for every $k$. Together with the lower bound $p_c(X)\geq p_m^-$, valid for every $X\geq1$ with mean $m$, this gives
\[
0<p_m^-\leq\inf\enstq{p_c(X)}{X\geq 1 \text{ a.s.},\ \Ec{X}=m}\leq\limsup_{k\to\infty}p_c(X_k^{(m)})\leq p_m^+.
\]
Hence the infimum is positive. Moreover, the asymptotic bounds on $p_m^-$ and $p_m^+$ show that both are equivalent to $\frac{2-m}{m+2}$ as $m\uparrow2$, and so is the infimum. Since $\frac{2-m}{m+2}\sim\frac{2-m}{4}$, this proves~\ref{it:pc without leaves m<2}.
\end{proof}

\subsection{A local perturbation lemma}
The counter-examples of Theorem~\ref{thm:contre exemple ordre convexe support} and Theorem~\ref{thm:contre-exemple m proche de 1} both rely on the following lemma.
\begin{lemma}\label{lem:local spread}
Let $X$ be a $\Z_+$-valued random variable with law $\mu$ and mean $m>1$, and let $k\geq1$ be such that $\mu(k)>0$. For $\epsilon\in\intervalleof{0}{\mu(k)/2}$, let $Y_\epsilon$ be a random variable with law
\begin{align*}
	\mu +  \epsilon \left(\delta_{k-1}-2\delta_k + \delta_{k+1}\right).
\end{align*}
Then $X\prec_{\mathrm{cx}}Y_\epsilon$. If moreover $p_c(X)^2\geq\frac{2}{k+1}$, then $p_c(X)<p_c(Y_\epsilon)$.
\end{lemma}
\begin{proof}
For every convex $f:\Z_+\to\R$, we have $\Ec{f(Y_\epsilon)}-\Ec{f(X)}=\epsilon\left(f(k+1)-2f(k)+f(k-1)\right)\geq0$, so that $X \preceq_{\mathrm{cx}}Y_\epsilon$. 
Since $\epsilon>0$, the laws of $X$ and $Y_\epsilon$ differ, hence $X\prec_{\mathrm{cx}}Y_\epsilon$.
Moreover, denoting by $\psi_\epsilon$ the generating function of $Y_\epsilon$, a simple computation gives $\psi_\epsilon(x) = \phi(x) + \epsilon x^{k-1}(x-1)^2$ so that the difference of their derivatives at any $x \in \intervalleoo{0}{1}$ can be written as
\begin{align}\label{eq:difference derivees contre-exemple}
	\psi_\epsilon'(x) - \phi'(x)=  \epsilon \cdot (1-x)x^{k-2}(k+1)\cdot \left(\frac{k-1}{k+1}-x\right).
\end{align}
Note that for $x\in\intervalleoo{0}{1}$ the sign of this difference is that of $\frac{k-1}{k+1}-x$.

Recall the spirit of the proof of Theorem~\ref{thm:convex order + bound support + large enough mean = comparison}: by Proposition~\ref{prop:criterion for comparison of critical points}, the ordering between the critical parameters $p_c(X)$ and $p_c(Y_\epsilon)$ is determined by the sign of $\psi_\epsilon'(f_{Y_\epsilon}(p_c(X))) - \phi'(f_X(p_c(X)))$.
As in~\eqref{eq:writing as good term plus bad term}, we write this difference as the sum of the good term $\psi_\epsilon'(f_{Y_\epsilon}(p_c(X))) - \phi'(f_{Y_\epsilon}(p_c(X)))$ and the bad term $\phi'(f_{Y_\epsilon}(p_c(X))) - \phi'(f_X(p_c(X)))$.
As before, from the fact that $X \preceq_{\mathrm{cx}}Y_\epsilon$, the bad term is non-negative: indeed by~\eqref{eq:inequality C p phi} we have $f_X(p) \leq f_{Y_\epsilon}(p)$ for all $p$ and the function $\phi'$ is non-decreasing.

Now, in the proof of Theorem~\ref{thm:convex order + bound support + large enough mean = comparison}, we crucially use the fact that the support of $X$ is bounded above by some value in order to show that the good term is negative enough to compensate the bad term.
On the contrary, here we are interested in showing that, under the appropriate assumptions, the good term is actually positive, so that the sum of the two terms is positive and $p_c(X)<p_c(Y_\epsilon)$.
Using \eqref{eq:difference derivees contre-exemple} we get
\begin{align*}
	\mathrm{sgn}\left(\psi_\epsilon'(f_{Y_\epsilon}(p_c(X))) - \phi'(f_{Y_\epsilon}(p_c(X)))\right) = \mathrm{sgn}\left(\frac{k-1}{k+1}- f_{Y_\epsilon}(p_c(X)) \right).
\end{align*}
It just remains to check that the assumption $p_c(X)^2\geq\frac{2}{k+1}$ is sufficient for $f_{Y_\epsilon}(p_c(X))<\frac{k-1}{k+1}$ to hold.
Recall that $p_c(X)\in\intervalleoo{0}{\frac1{m^2}}$ by Theorem~\ref{prop:the critical percolation parameter is the solution to an equation}.
Hence Lemma~\ref{lem:proprietes f_X et F_m}\ref{it:universal bound Fxp} applied to $Y_\epsilon$ gives $f_{Y_\epsilon}(p_c(X))<1-p_c(X)^2\leq1-\frac{2}{k+1}=\frac{k-1}{k+1}$, which finishes the proof.
\end{proof}

\subsection{Proof of Theorem~\ref{thm:contre exemple ordre convexe support} and Theorem~\ref{thm:contre-exemple m proche de 1}}
\begin{proof}[Proof of Theorem~\ref{thm:contre exemple ordre convexe support}]
In both parts, the random variable $Y$ is one of the random variables $Y_\epsilon$ of Lemma~\ref{lem:local spread}, for a well-chosen integer $k$.
Indeed, for a $\Z_+$-valued random variable $X$ with mean $m$ and an integer $k\geq1$, this lemma gives $X\prec_{\mathrm{cx}}Y_\epsilon$ and $p_c(X)<p_c(Y_\epsilon)$ for every $\epsilon\in\intervalleof{0}{\Pp{X=k}/2}$ as soon as
\begin{equation}\label{eq:condition local spread}
\Pp{X=k}>0\qquad\text{and}\qquad p_c(X)^2\geq\frac{2}{k+1}.
\end{equation}

We start with point~\ref{it:counter-example unbounded support}. 
Since $p_c(X)>0$ and the support of $X$ is unbounded, there exists $k\geq1$ in the support of $X$ such that $p_c(X)^2\geq\frac{2}{k+1}$, so that~\eqref{eq:condition local spread} holds.

We now prove point~\ref{it:counter-example bounded support}. 
We construct a random variable $X$ with mean $m$ such that~\eqref{eq:condition local spread} holds for an explicit integer $k=N'(m)$, which is the largest element of the support of $X$.
In Step~\ref{step:lower bound pc}, we prove a lower bound on $p_c(X)$ which holds for every $X$ with mean $m$ satisfying a given upper bound on $\Ec{X(X-1)}$, which in particular holds for $X^{(m)}$.
In Step~\ref{step:construction X}, we define $N'(m)$ so that the second condition in~\eqref{eq:condition local spread} holds with $k=N'(m)$ as soon as $X$ satisfies this moment bound.
We then construct $X$ as a mixture of the laws of $X^{(m)}$ and $X^{(m)}_{N'(m)}$, which both have mean $m$, with a small weight on the latter: the first condition in~\eqref{eq:condition local spread} holds since $X^{(m)}_{N'(m)}$ charges $N'(m)$, and the small weight ensures that $X$ satisfies the moment bound.
We conclude in Step~\ref{step:conclusion}, and Step~\ref{step:asymptotics N'} is devoted to the asymptotics of $N'(m)$.

\step{a lower bound on $p_c$}\label{step:lower bound pc}
Suppose $X$ is a $\Z_+$-valued random variable with mean $m$ and finite second moment, and denote by $\phi$ its generating function. 
The function $\phi'$ is convex on $\intervalleff{0}{1}$ and $\phi''(1)=\Ec{X(X-1)}$, so that, for $p\in\intervalleoo{0}{\frac1m}$,
\[
\phi'(f_X(p))\geq m-\Ec{X(X-1)}\left(1-f_X(p)\right).
\]
Now, we impose the following extra condition on $X$
\begin{equation}\label{hyp:moment}
\Ec{X(X-1)}\leq\lceil m\rceil(\lceil m\rceil-1)+1,
\end{equation}
which entails, using Lemma~\ref{lem:proprietes f_X et F_m}\ref{it:universal bound Fxp},
\begin{align*}
\phi'(f_X(p)) \geq m-\frac{\left(\lceil m\rceil(\lceil m\rceil-1)+1\right)\,p^2}{1-p(1-p)m}.
\end{align*}
We then set, for every $p\in\intervallefo{0}{\frac1m}$,
\begin{equation}\label{eq:definition vm}
v_m(p):=\max\left(0,\,m-\frac{\left(\lceil m\rceil(\lceil m\rceil-1)+1\right)\,p^2}{1-p(1-p)m}\right),
\end{equation}
so that, since $\phi'$ is non-negative, we have $\phi'(f_X(p))\geq v_m(p)$. 

The function $F_m(p,\cdot)$ is increasing on $\intervallefo{0}{\infty}$ by Lemma~\ref{lem:proprietes f_X et F_m}\ref{it:p,v mapsto Fpv is continuous}, so that $F_m\left(p,\phi'(f_X(p))\right)\geq F_m\left(p,v_m(p)\right)$. 
By Lemma~\ref{lem:p is supercritical iff F<0}, there is therefore no percolation at $p$ as soon as $F_m(p,v_m(p))>0$. 
This leads us to define
\[
\underline{p}(m):=\inf\enstq{p\in\intervalleoo{0}{\tfrac1m}}{F_m(p,v_m(p))\leq0}.
\]
The set in this definition is non-empty: since $v_m\left(\frac1{m^2}\right)\leq m$ and $F_m$ is increasing in its second variable, we have
\[
F_m\left(\frac1{m^2},v_m\left(\frac1{m^2}\right)\right)\leq F_m\left(\frac1{m^2},m\right)=(m-1)\left[\left(1-\frac1{m^2}\right)\left(1-\frac1{m^4}\right)-1\right]<0,
\]
so that it contains $\frac1{m^2}$.
Hence every $p\in\intervalleoo{0}{\underline{p}(m)}$ satisfies $F_m(p,v_m(p))>0$, so that there is no percolation at $p$. 
This shows that $p_c(X)\geq\underline{p}(m)$.

\step{construction of $X$}\label{step:construction X}
The function $p\mapsto F_m(p,v_m(p))$ is continuous on $\intervallefo{0}{\frac1m}$ by Lemma~\ref{lem:proprietes f_X et F_m}\ref{it:p,v mapsto Fpv is continuous}, and its value at $p=0$ is $F_m(0,v_m(0))=F_m(0,m)=1$, so that $\underline{p}(m)>0$. 
We can therefore define
\begin{equation}\label{eq:definition N'}
N'(m):=\left\lceil\frac{2}{\underline{p}(m)^2}\right\rceil ,
\end{equation}
which is an integer larger than $m$: since $\underline{p}(m)<\frac1m$, we have $N'(m)>2m^2>m$. 
Recall from~\eqref{eq:def loi X^{(m)}} that the law of $X^{(m)}$ is $(1-\{m\})\delta_{\lfloor m\rfloor}+\{m\}\delta_{\lfloor m\rfloor+1}$, and from the beginning of this section that the law of $X_{N'(m)}^{(m)}$ is $\frac{N'(m)-m}{N'(m)-1}\delta_1+\frac{m-1}{N'(m)-1}\delta_{N'(m)}$. 
For $\eta\in\intervalleoo{0}{1}$, let $X$ be a random variable whose law is the mixture
\[
(1-\eta)\left((1-\{m\})\delta_{\lfloor m\rfloor}+\{m\}\delta_{\lfloor m\rfloor+1}\right)+\eta\left(\frac{N'(m)-m}{N'(m)-1}\delta_1+\frac{m-1}{N'(m)-1}\delta_{N'(m)}\right).
\]
As a mixture of two laws with mean $m$, the law of $X$ has mean $m$. 
Moreover, since $N'(m)$ is an integer larger than $m$, we have $\supp X\subset\{1,\dots,N'(m)\}$ and $\Pp{X=N'(m)}\geq\eta\frac{m-1}{N'(m)-1}>0$. 
Finally, the quantity
\[
\Ec{X(X-1)}=(1-\eta)\,\Ec{X^{(m)}\left(X^{(m)}-1\right)}+\eta\,\Ec{X_{N'(m)}^{(m)}\left(X_{N'(m)}^{(m)}-1\right)}
\]
tends to $\Ec{X^{(m)}(X^{(m)}-1)}$ as $\eta\to0$, and this limit is at most $\lceil m\rceil(\lceil m\rceil-1)$ since $X^{(m)}\leq\lceil m\rceil$. 
Hence $X$ satisfies~\eqref{hyp:moment} for $\eta$ small enough, and we fix such an $\eta$ from now on.

\step{conclusion}\label{step:conclusion}
By Step~\ref{step:lower bound pc} and the definition~\eqref{eq:definition N'} of $N'(m)$, we have
\[
p_c(X)^2\geq\underline{p}(m)^2\geq\frac{2}{N'(m)}\geq\frac{2}{N'(m)+1},
\]
and $\Pp{X=N'(m)}>0$ by Step~\ref{step:construction X}. 
Hence~\eqref{eq:condition local spread} holds with $k=N'(m)$.

\step{asymptotics of $N'(m)$}\label{step:asymptotics N'}
In this step only, we assume that $m>3$. 
Let $\theta\in\intervalleof{0}{1-\frac3m}$. 
We have $\frac{\theta}{m^2}\left(1-\frac{\theta}{m^2}\right)m\leq\frac1m\leq\frac13$, so that the denominator in the definition~\eqref{eq:definition vm} of $v_m\left(\frac{\theta}{m^2}\right)$ is at least $\frac23$. 
Since $\lceil m\rceil(\lceil m\rceil-1)+1\leq m^2+m+1\leq\frac{13}{9}m^2$, we get
\[
v_m\left(\frac{\theta}{m^2}\right)\geq m-\frac{13}{9}m^2\cdot\frac{\theta^2}{m^4}\cdot\frac32=m-\frac{13\theta^2}{6m^2}.
\]
For $p\in\intervalleff{0}{\frac1m}$ and $v\in\intervalleff{0}{m}$, we have $(1-p)^2\left(1-mp+p^2(m-1)v\right)\geq(1-2p)(1-mp)\geq1-(m+2)p$, hence $F_m(p,v)\geq v-(m+2)p\,m-(m-1)$. 
Taking $p=\frac{\theta}{m^2}$ and $v=v_m\left(\frac{\theta}{m^2}\right)$, we get
\[
F_m\left(\frac{\theta}{m^2},v_m\left(\frac{\theta}{m^2}\right)\right)\geq v_m\left(\frac{\theta}{m^2}\right)-\left(1+\frac2m\right)\theta-(m-1)\geq 1-\theta-\frac{2\theta}{m}-\frac{13\theta^2}{6m^2}.
\]
The right-hand side is decreasing in $\theta$, and for $\theta=1-\frac3m$ it is at least $\frac1m-\frac{13}{6m^2}>0$. 
Hence $F_m(p,v_m(p))>0$ for every $p\in\intervalleof{0}{\frac{1-3/m}{m^2}}$, so that $\underline{p}(m)\geq\frac{1-3/m}{m^2}$ and
\[
N'(m)\leq\frac{2m^4}{(1-3/m)^2}+1=(2+o(1))\,m^4\qquad\text{as }m\to\infty.
\]
This is the bound on $N'(m)$ stated in the theorem, which finishes the proof.
\end{proof}

\begin{remark}\label{rem:N'(2) numerics}
For $m=2$ we find $\underline{p}(2)\approx0.1549$ and $N'(2)=84$.
\end{remark}

\begin{proof}[Proof of Theorem~\ref{thm:contre-exemple m proche de 1}]
For $m\in\intervalleff{1}{1.1}$ and $p\in\intervalleff{0}{0.9}$, since $1-p(1-p)m-p^2=(1-p)(1+p-pm)$, the quantity $1-\frac{p^2}{1-p(1-p)m}$ is positive. 
Hence the right-hand side of~\eqref{eq:lower bound F without leaves}, which makes sense for $m=1$ as well, is a continuous function of $(m,p)$ on the compact set $\intervalleff{1}{1.1}\times\intervalleff{0}{0.9}$. 
For $m=1$, it is equal to $F_1(p,1)=(1-p)^3\geq 0.1^3>0$. By uniform continuity, there exists $\delta\in\intervalleoo{0}{0.1}$ such that it is positive for every $m\in\intervallefo{1}{1+\delta}$ and every $p\in\intervalleff{0}{0.9}$. Note that $0.9<\frac1m$ for such $m$.

Let $m\in\intervalleoo{1}{1+\delta}$. 
The random variable $X^{(m)}$ satisfies $X^{(m)}\geq1$ and $\Pp{X^{(m)}=2}=m-1>0$. Denoting by $\phi$ its generating function, Lemma~\ref{lem:lower bound without leaves} gives $F_m\left(p,\phi'(f_{X^{(m)}}(p))\right)>0$ for every $p\in\intervalleof{0}{0.9}$. 
By Lemma~\ref{lem:p is supercritical iff F<0}, we get $p_c(X^{(m)})\geq 0.9$, so that $p_c(X^{(m)})^2\geq 0.81>\frac23$. 
Lemma~\ref{lem:local spread}, applied to $X^{(m)}$ with $k=2$, for which $\frac{2}{k+1}=\frac23$, gives the result. 
(The same argument applies to every $X\geq 1$ with mean $m$ and $\Pp{X=2}>0$.)
\end{proof}
\appendix
\section{Annealed and quenched critical parameters}\label{sec:annealed quenched}
Write
\begin{align*}
	\theta(p)= \Pp{\emptyset \overset{p}{\leftrightarrow} \infty}
\end{align*}
for the probability that the root belongs to an infinite $p$-cluster on $\GF(T)$, where, as in the introduction, $\GF(T)$ does not contain the vertex $\hat\emptyset$. 
By the discussion at the end of Section~\ref{sec:notation}, we have $\theta(p)>0$ if and only if $\Pp{|\mathcal C_p|=\infty}>0$.
We can also introduce the quantity
\begin{align*}
	\theta(p,T)= \Ppsq{\emptyset \overset{p}{\leftrightarrow} \infty}{T},
\end{align*}
which is the same quantity but this time \emph{conditionally on the tree $T$}.
\begin{lemma}\label{lem:annealed quenched}
For any $p\in \intervalleff{0}{1}$, if $\theta(p)>0$ then almost surely on $\{|T| = \infty\}$ we have
	\begin{align*}
		\theta(p,T)>0.
	\end{align*}
\end{lemma}
\begin{proof}
Let $p\in \intervalleff{0}{1}$ and assume that $\theta(p)>0$.
First note that we have the inclusion of events
\begin{align}\label{eq:inclusion of events T is infinite}
	\{\theta(p,T)>0\} \subset \{|T|=\infty\}.
\end{align}
In order to prove the lemma, we just need to prove that those two events have the same probability.
Denote by $\alpha(p)$ the probability $\Pp{\theta(p,T)>0}$ of the event on the left-hand side.
By our assumption on $p$, we have $\Ec{\theta(p,T)}= \theta(p)>0$, so that we have $\alpha(p) >0$.

We introduce a bit of notation: for any integer $n \geq0$, we let $\Xi_n$ denote the set of vertices of $T$ at height $n$.
Recall that, for $v\in \Xi_n$, we denote by $T_v$ the subtree of the descendants of $v$, which we see as a tree rooted at $v$.
Note that for any such $v$, in order for the root $\emptyset$ to be connected to $\infty$, it is sufficient that, at the same time
\begin{itemize}
	\item all the short edges on the path from $\emptyset$ to $v$ are open,
	\item the vertex $v$ is connected to $\infty$ in $\GF(T_v)$, that is using only edges with both endpoints in $T_v$.
\end{itemize}
Conditionally on $T$, these two events are independent, since they depend on disjoint sets of edges, so we may conclude that
\begin{align*}
	\theta(p,T) \geq p^n \theta(p, T_v),
\end{align*}
so that $\theta(p,T) $ can only vanish when all the $(\theta(p, T_v))_{v \in \Xi_n}$ vanish as well.
Now, conditionally on $\Xi_n$, the trees $(T_v)_{v \in \Xi_n}$ are independent Bienaymé trees with the same law as $T$.
This entails that
\begin{align}\label{eq:decomposition probability at height n}
	\alpha(p)=\Pp{\theta(p,T)>0} \geq \Ec{ 1 -\Ppsq{\forall v \in \Xi_n, \ \theta(p,T_v) =0}{\Xi_n} } = \Ec{1 - (1-\alpha(p))^{|\Xi_n|}}.
\end{align}
Now, since $T$ is a (non-degenerate) Bienaymé tree,
its generation sizes $(|\Xi_n|)_{n\geq 0}$ form a Markov chain on $\Z_+$ that can only be absorbed at $0$ in finite time or tend to $\infty$.
Hence  $|\Xi_n| \rightarrow \infty$ almost surely as $n\rightarrow \infty$ on $\{|T|=\infty\}$ and $|\Xi_n|=0$ for $n$ large enough on $\{|T|<\infty\}$.
Since $\alpha(p)>0$ we get that almost surely
\begin{align*}
	 1 - (1-\alpha(p))^{|\Xi_n|} \underset{n\rightarrow \infty}{\rightarrow} 1-  \ind{|T|< \infty} =  \ind{|T|= \infty}.
\end{align*}
Since the random variables above are bounded, we can take expectations in the limit above.
Recalling \eqref{eq:decomposition probability at height n}, this proves that
\begin{align*}
\alpha(p)	=\Pp{\theta(p,T)>0} \geq  \Pp{|T|=\infty},
\end{align*}
which, combined with \eqref{eq:inclusion of events T is infinite}, implies the statement of the lemma.
\end{proof}
Let us deduce the claim made in the introduction, namely that, almost surely on $\{|T|=\infty\}$,
\[
\inf\enstq{p\in\intervalleff{0}{1}}{\theta(p,T)>0}=p_c(X).
\]
By definition of $p_c(X)$, and since $p\mapsto\theta(p)$ is non-decreasing, we have $\theta(p)>0$ for $p>p_c(X)$ and $\theta(p)=0$ for $p<p_c(X)$. 
For every $p>p_c(X)$, Lemma~\ref{lem:annealed quenched} gives $\theta(p,T)>0$ almost surely on $\{|T|=\infty\}$. 
For every $p<p_c(X)$, we have $\Ec{\theta(p,T)}=\theta(p)=0$, so that $\theta(p,T)=0$ almost surely. 
Applying these two facts to all rational $p$, and using that $p\mapsto\theta(p,T)$ is non-decreasing, we get the claim.

\end{document}